\RequirePackage[l2tabu, orthodox]{nag}

\documentclass[12pt,reqno]{amsart}
\usepackage{fullpage,url,amssymb,enumerate,colonequals,graphicx}
\usepackage{caption}
\usepackage{subcaption}
\usepackage{multirow}
\usepackage{float}
\usepackage{mathrsfs} % for \mathscr (script letters)

\usepackage[dvipsnames,xcdraw,hyperref]{xcolor}

\newcommand{\defi}[1]{\textsf{#1}} % for defined terms

\newcommand{\C}{\mathbb{C}}
\newcommand{\F}{\mathbb{F}}

\newcommand{\PP}{\mathbb{P}}
\newcommand{\Q}{\mathbb{Q}}
\newcommand{\Qtilde}{\tilde{Q}}
\newcommand{\R}{\mathbb{R}}
\newcommand{\Z}{\mathbb{Z}}
\newcommand{\Qbar}{{\overline{\Q}}}

\newcommand{\Fbar}{{\overline{\F}}}

\newcommand{\kk}{k}
\newcommand{\kap}{\kappa}

\newcommand{\diff}{\textup{diff}}
\newcommand{\final}{\textup{final}}

\newcommand{\Ia}{I_{\textup{attained}}}
\newcommand{\Itrue}{I_{\textup{true}}}
\newcommand{\Iinner}{I_{\textup{inner}}}
\newcommand{\Iouter}{I_{\textup{outer}}}
\newcommand{\Isimple}{I_{\textup{simple}}}
\newcommand{\IWeil}{I_{\textup{Weil}}}

\newcommand{\calA}{\mathcal{A}}

\newcommand{\calC}{\mathcal{C}}

\newcommand{\calG}{\mathcal{G}}

\newcommand{\calJ}{\mathcal{J}}

\newcommand{\calM}{\mathcal{M}}

\newcommand{\calS}{\mathcal{S}}

\newcommand{\FF}{\mathscr{F}}

\newcommand{\MM}{\mathscr{M}}

\DeclareMathOperator{\Gal}{Gal}

\DeclareMathOperator{\ord}{ord}

\DeclareMathOperator{\re}{Re}

\DeclareMathOperator{\supp}{supp}

\newcommand{\del}{\partial}
\newcommand{\surjects}{\twoheadrightarrow}
\newcommand{\To}{\longrightarrow}
\newcommand{\union}{\cup} % binary union
\newtheorem{theorem}{Theorem}[section]
\newtheorem{lemma}[theorem]{Lemma}
\newtheorem{corollary}[theorem]{Corollary}
\newtheorem{proposition}[theorem]{Proposition}

\theoremstyle{definition}
\newtheorem{definition}[theorem]{Definition}

\newtheorem{construction}[theorem]{Construction}

\theoremstyle{remark}
\newtheorem{remark}[theorem]{Remark}

\usepackage{scalerel,stackengine}
\stackMath
\newcommand\reallywidehat[1]{%
\savestack{\tmpbox}{\stretchto{%
  \scaleto{%
    \scalerel*[\widthof{\ensuremath{#1}}]{\kern.1pt\mathchar"0362\kern.1pt}%
    {\rule{0ex}{\textheight}}%WIDTH-LIMITED CIRCUMFLEX
  }{\textheight}% 
}{2.4ex}}%
\stackon[-8.5pt]{#1}{\tmpbox}%
}
\makeatletter
\g@addto@macro\bfseries{\boldmath} % This makes math in section titles bold.
\makeatother

\usepackage{tikz}
\usepackage{xcolor}
\usepackage{color}
\definecolor{mylinkcolor}{rgb}{0.5,0.0,0.0}
\definecolor{myurlcolor}{rgb}{0.0,0.0,0.75}
\usepackage[
	colorlinks, urlcolor=myurlcolor,citecolor=myurlcolor,linkcolor=mylinkcolor,
	pdfauthor={Raymond van Bommel, Edgar Costa, Wanlin Li, Bjorn Poonen, Alexander Smith}, % add other authors
	pdftitle={Abelian varieties of prescribed order over finite fields},
]{hyperref}
\usepackage{csquotes} % removes one warning of biblatex
\usepackage[style=ieee-alphabetic,backend=biber, sorting=nyt, doi=false,isbn=false,url=false]{biblatex}
\newbibmacro{string+doiurlisbn}[1]{%
  \iffieldundef{doi}{%
    \iffieldundef{url}{%
      \iffieldundef{isbn}{%
        \iffieldundef{issn}{%
          #1%
        }{%
          \href{http://books.google.com/books?vid=ISSN\thefield{issn}}{#1}%
        }%
      }{%
        \href{http://books.google.com/books?vid=ISBN\thefield{isbn}}{#1}%
      }%
    }{%
      \href{\thefield{url}}{#1}%
    }%
  }{%
    \href{http://dx.doi.org/\thefield{doi}}{#1}%
  }%
}

\DeclareFieldFormat{title}{\usebibmacro{string+doiurlisbn}{\mkbibemph{#1}}}
\DeclareFieldFormat[article,incollection,inproceedings]{title}%
    {\usebibmacro{string+doiurlisbn}{\mkbibquote{#1}}}
\usepackage{microtype}  % This adjusts spacing between words so as to improve the probability of having line breaks in good places.

\makeatletter
\@namedef{subjclassname@2020}{%
  \textup{2020} Mathematics Subject Classification}
\makeatother

\begin{document}

\title[short paper title]{Abelian varieties of prescribed order over finite fields}
% 11G10 Abelian varieties of dimension > 1
% 11G25 Varieties over finite and local fields
% 11Y99 Computational number theory
% 14G15 Finite ground fields in algebraic geometry
% 14K15 Arithmetic ground fields for abelian varieties
% 31A15 potential theory

\author{Raymond van Bommel}
\address{Department of Mathematics, Massachusetts Institute of Technology, Cambridge, MA 02139-4307, USA}
\email{bommel@mit.edu}
\urladdr{\url{https://www.raymondvanbommel.nl/}}

\author{Edgar Costa}
\address{Department of Mathematics, Massachusetts Institute of Technology, Cambridge, MA 02139-4307, USA}
\email{edgarc@mit.edu}
\urladdr{\url{https://edgarcosta.org}}

\author{Wanlin Li}
\address{Centre de Recherches Math\'ematiques,
         Universit\'e de Montr\'eal, 2920 Chemin de la tour,
         Montr\'eal (Qu\'ebec) H3T 1J4, Canada}
\email{liwanlin@crm.umontreal.ca}
\urladdr{\url{http://www.crm.umontreal.ca/~liwanlin/}}

\author{Bjorn Poonen}
\address{Department of Mathematics, Massachusetts Institute of Technology, Cambridge, MA 02139-4307, USA}
\email{poonen@math.mit.edu}
\urladdr{\url{https://math.mit.edu/~poonen/}}

\author{Alexander Smith}
\address{Department of Mathematics, Massachusetts Institute of Technology, Cambridge, MA 02139-4307, USA}
\email{adsmi@mit.edu}

\thanks{R.B., E.C., and B.P.\ were supported by Simons Foundation grant \#550033.
B.P.\ was supported also in part by National Science Foundation grant DMS-1601946 and Simons Foundation grant \#402472. A.S.\ was supported in part by National Science Foundation grant DMS-2002011.  W.L.\ is supported by a fellowship from the Centre de recherches math\'ematiques (CRM) and the Institut des sciences mathématiques (ISM), and funds from the Natural Sciences and Engineering Research Council of Canada (NSERC) and the Fonds de recherche du Qu\'ebec - Nature et technologies (FRQNT)}

\date{June 24, 2021}

\begin{abstract}
Given a prime power $q$ and $n \gg 1$, we prove that every integer in a large subinterval of the Hasse--Weil interval $[(\sqrt{q}-1)^{2n},(\sqrt{q}+1)^{2n}]$ is $\#A(\F_q)$ for some geometrically simple ordinary principally polarized abelian variety $A$ of dimension $n$ over $\F_q$.  
As a consequence, we generalize a result of Howe and Kedlaya for $\F_2$ to show that for each prime power $q$, every sufficiently large positive integer is realizable, i.e., $\#A(\F_q)$ for some abelian variety $A$ over $\F_q$.
Our result also improves upon the best known constructions of sequences of simple abelian varieties with point counts towards the extremes of the Hasse--Weil interval.
A separate argument determines, for fixed $n$, the largest subinterval of the Hasse--Weil interval consisting of realizable integers, asymptotically as $q \to \infty$; this gives an asymptotically optimal improvement of a 1998
theorem of DiPippo and Howe.
Our methods are effective: We prove that if $q \le 5$, then every positive integer is realizable,
and for arbitrary $q$, every positive integer $\ge q^{3 \sqrt{q} \log q}$ is realizable.
\end{abstract}

\maketitle

%****************************************************************************
\section{Introduction}\label{S:introduction}

\subsection{Orders of abelian varieties over a finite field}

By work of Weil (a consequence of \cite[pp.~70--71]{Weil-1948-Sur} and \cite[pp.~137--138]{Weil-1948-Varietes}, generalizing
\cite[p.~206]{Hasse-1936c}), 
if $A$ is an abelian variety of dimension $n$ over a finite field $\F_q$, then $\#A(\F_q)$ lies in the interval
\begin{equation}
\label{E:Weil bound}
  \bigl[ \; \left( q - 2q^{1/2} + 1 \right)^n \;,\; \left( q + 2q^{1/2} + 1 \right)^n \; \bigr].
\end{equation}
We prove an almost-converse (compare \eqref{E:Weil bound} and \eqref{E:simplified interval}):
\begin{theorem}
\label{T:almost-converse}
Fix a prime power $q$.
Let $\tau(x) = x+\sqrt{x^2-1}$.
Let $I$ be a closed interval contained in
\begin{equation}\label{E:I attainable}
\Ia \colonequals \bigl( \; \tau(q/2 - q^{1/2} + 3/2) \; , \; \tau(q/2 + q^{1/2} - 1/2)  \; \bigr). \\
\end{equation}
For $n$ sufficiently large, 
if $m$ is a positive integer with $m^{1/n} \in I$,
then there exists an $n$-dimensional abelian variety $A$ with $\#A(\F_q)=m$.
Moreover, $A$ can be chosen to be ordinary, geometrically simple, and principally polarized.
\end{theorem}

\begin{corollary}
\label{C:simplified interval}
Fix a prime power $q$. Then, for $n$ sufficiently large, every integer in the interval
\begin{equation}
\label{E:simplified interval}
    \bigl[ \; \left(q - 2q^{1/2} + 3 - q^{-1} \right)^n \; , \; \left(q + 2q^{1/2} - 1 - q^{-1}\right)^n  \; \bigr]
\end{equation}
is $\#A(\F_q)$ for some geometrically simple ordinary principally polarized abelian variety $A$ of dimension $n$ over $\F_q$.
\end{corollary}

The interval \eqref{E:simplified interval} in Corollary~\ref{C:simplified interval} contains $[q^n,q^{n+1}]$ if $n$ is large enough,
so Corollary~\ref{C:simplified interval} implies the following:

\begin{corollary}
\label{C:HK for q}
Fix a prime power $q$.
Every sufficiently large positive integer is $\#A(\F_q)$
for some geometrically simple ordinary principally polarized abelian variety $A$ over $\F_q$.
\end{corollary}

Corollary~\ref{C:HK for q} answers a question of Howe and Kedlaya, who proved that every positive integer
is the order of an ordinary abelian variety over $\F_2$ \cite[Theorem~1]{Howe-Kedlaya-preprint}.
For effective versions, see Section~\ref{S:introduction effective}.

\begin{remark}
Marseglia and Springer refined \cite{Howe-Kedlaya-preprint}
to prove that every finite abelian group is isomorphic to $A(\F_2)$ for some ordinary abelian variety $A$ over $\F_2$ \cite{Marseglia-Springer-preprint}. Our Corollary~\ref{C:HK for q} combined with~\cite[Theorem~4.2]{Marseglia-Springer-preprint} implies that for any fixed $q$, every cyclic group of sufficiently large order is isomorphic to $A(\F_q)$ for some ordinary abelian variety $A$ over $\F_q$.
\end{remark}

Throughout, $p$ denotes the characteristic of $\F_q$.

\begin{remark}
Theorem~\ref{T:almost-converse} can be extended to produce non-ordinary abelian varieties.
First, define the \defi{$p$-rank} of an $n$-dimensional abelian variety $A$ over $\F_q$ to be  
the integer $\dim_{\F_p} A[p](\Fbar_q)$ in $[0,n]$.
For example, $A$ is ordinary if and only if the $p$-rank is $n$.
Then Theorem~\ref{T:almost-converse} holds with ``ordinary'' replaced by ``of prescribed $p$-rank $r$'' for
any $r \in [0,n]$, provided that when $r=0$, we assume $m \equiv 1 \pmod{p}$; see Remark~\ref{R:congruences for p-rank}.
\end{remark}

\begin{remark}
It may be that Theorem~\ref{T:almost-converse} holds for an interval larger than $\Ia$.
There is a largest open interval $\Itrue$ containing $q$ for which Theorem~\ref{T:almost-converse} holds.
\end{remark}

\subsection{Extreme point counts for simple abelian varieties}
\label{S:extreme}

Other authors have studied the extreme values of $\#A(\F_q)^{1/\dim A}$ without trying to realize
every order in between.
Following~\cite{Kadets-2021}, let $\calA_q$ be the set of simple abelian varieties over $\F_q$ up to isogeny and consider
\[
    \Isimple \colonequals \bigl[ \;  \liminf_{A \in \calA_q} \#A(\F_q)^{1/\dim A} \;,\; \limsup_{A \in \calA_q} \#A(\F_q)^{1/\dim A} \; \bigr].
\]
(If one did not require simplicity and take $\limsup$ and $\liminf$, then for square $q$ the minimum and maximum would be achieved by elliptic curves of order $q \pm 2q^{1/2} + 1$ and their powers.)
Then
\[   
   \Ia \subseteq \Itrue \subseteq \Isimple \subseteq \IWeil \colonequals [\; q-2q^{1/2}+1 \;,\; q+2q^{1/2}+1\;].
\]
Aubry, Haloui and Lachaud \cite[Corollaries 2.2 and~2.14]{Aubry-Haloui-Lachaud-2013} and Kadets~\cite[Theorem~1.8]{Kadets-2021} 
found inner and outer bounds $\Iinner,\Iouter$ for $\Isimple$: 
\begin{equation}
\label{E:inner and outer}
    \Bigl[ \; 
    q - \lfloor 2q^{1/2} \rfloor + 3 % = \lceil (q^{1/2}-1)^2 \rceil + 2 
    \;,\; 
    q + \lfloor 2q^{1/2} \rfloor - 1 - q^{-1} %= \lfloor (q^{1/2}+1)^2 \rfloor - 2 - q^{-1} 
    \;\Bigr] 
    \; \subseteq \; \Isimple \; \subseteq \;
    \Bigl[\; 
     q - \lceil 2q^{1/2} \rceil + 2   % = \lfloor (q^{1/2}-1)^2 \rfloor + 1 
    \;,\; 
   q + \lceil 2q^{1/2} \rceil % = \lceil (q^{1/2}+1)^2 \rceil - 1 
    \; \Bigr].
\end{equation}
Our inner bound $\Ia$ for $\Isimple$ improves upon $\Iinner$,
but careful consideration shows that Kadets's argument yields
a better result than he claimed, an inner bound matching our $\Ia$
when $q$ is a square.

The following diagram shows $\Ia \subset \Iouter \subset \IWeil$, bounded by open dots, solid dots, and vertical bars, respectively. 
The endpoints of $\Itrue$ and $\Isimple$ are unknown, but they lie somewhere in the (closed) dashed intervals.

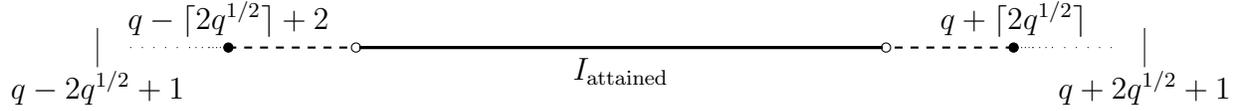
\begin{center}
\begin{tikzpicture}[scale = 0.87]
\draw[very thick] (-4,0) -- (4,0);
\draw (-8,-0.3) -- (-8,0.3);
\draw (8,-0.3) -- (8,0.3);
\draw [thick,dashed] (-6, 0) -- (-4.1,0);
\draw [thick,dashed] (6, 0) -- (4.1,0);

\draw [loosely dotted] (-7.5, 0) -- (-6.7,0);
\draw [dotted] (-6.6, 0) -- (-6.3,0);
\draw [densely dotted] (-6.3, 0) -- (-6,0);

\draw [loosely dotted] (7.5, 0) -- (6.7,0);
\draw [dotted] (6.6, 0) -- (6.3,0);
\draw [densely dotted] (6.3, 0) -- (6,0);

\node[below] at (-8,-0.2) {$q-2q^{1/2}+1$};
\node[below] at (8,-0.2) {$q+2q^{1/2}+1$};
\draw [fill] (-6,0) circle [radius=0.07];
\node[above] at (-6,0) {$q - \lceil 2q^{1/2} \rceil + 2$};
\draw [fill] (6,0) circle [radius=0.07];
\node[above] at (6,0) {$q + \lceil 2q^{1/2}\rceil$};

\draw (-4.05,0) circle [radius=0.07];
\draw  (4.05,0) circle [radius=0.07];
\node[below] at (0,0) {$\Ia$};

\end{tikzpicture}
\end{center}

\subsection{Strategy of proof}
\label{S:strategy}

Given an abelian variety $A$ over the finite field $\F_q$,
let $f_A(x) \in \Z[x]$ be the characteristic polynomial of the $q$-power Frobenius
acting on a Tate module $T_\ell A$.
Then $\#A(\F_q)=f_A(1)$.
Honda--Tate theory implies that for $f \in \Z[x]$,
we have $f=f_A$ for some ordinary $n$-dimensional abelian variety $A$ over $\F_q$
if and only if $f$ is monic of degree $2n$
with complex roots $\alpha_1,\bar{\alpha}_1,\ldots,\alpha_n,\bar{\alpha}_n$ satisfying $|\alpha_i|=q^{1/2}$,
and $p$ does not divide the coefficient of $x^n$.
Therefore, as in \cite{Howe-Kedlaya-preprint}, we need to find a polynomial $f$ satisfying these
conditions with a prescribed value of $f(1)$.

One ingredient that lets us go beyond \cite{Howe-Kedlaya-preprint}
is a lemma more general than \cite[Lemma~3.3.1]{DiPippo-Howe-1998} 
for constructing polynomials whose roots lie on the circle $|z|=q^{1/2}$ 
(Lemma~\ref{L:h to f}).
Using this lemma alone, we can give a quick proof of Corollary~\ref{C:HK for q}, 
if we omit ``geometrically simple'' and ``principally polarized'': see Section~\ref{S:simple proof}.

To force $A$ to be geometrically simple and principally polarized, 
we prove that it suffices to impose certain congruence 
conditions on the coefficients of $f$ (Proposition~\ref{P:good polynomial mod L}); 
unlike \cite[Lemma~3.3.1]{DiPippo-Howe-1998}, 
our Lemma~\ref{L:h to f} is robust enough to permit a wide enough range of values $f(1)$
even when such congruence conditions are imposed.
To prove Theorem~\ref{T:almost-converse}, 
we start with rescaled Chebyshev polynomials similar to those in~\cite{Kadets-2021} (Proposition~\ref{P:construction of P}),
but we improve on \cite{Kadets-2021} by temporarily allowing \emph{non-integral} real coefficients,
and later making adjustments to make the coefficients integral while preserving $f(1)$
and the bounds needed to apply Lemma~\ref{L:h to f}.
To obtain the widest interval of realizable values, 
we must adjust differently in three different ranges of exponents, 
and the adjustments do something more elaborate than changing one coefficient at a time; see Section~\ref{sec:range}.

Although we do not know if the bounds in Theorem~\ref{T:almost-converse} are sharp, 
Appendix~\ref{A:potential theory} shows that the rescaled Chebyshev polynomials
are asymptotically optimal for our \emph{method}.

\subsection{Large \texorpdfstring{$q$}{q} limit}

So far we have discussed the possibilities for $\#A(\F_q)$
for an $n$-dimensional abelian variety over a fixed finite field $\F_q$, as $n \to \infty$.
We also obtain a sharp asymptotic for the possibilities for fixed $n$ as $q \to \infty$:

\begin{theorem}
\label{T:large q}
Fix $n \ge 3$.
Let $\lambda_1 = 2n - \sqrt{\tfrac{2n}{n-1}}$.
Then the largest interval in which every integer is $\#A(\F_q)$
for some $n$-dimensional abelian variety $A$ over $\F_q$ has the form
\begin{equation}
\label{E:interval for large q}
   \Bigl[\; q^n - \lambda_1 q^{n-1/2} + o(q^{n-1/2}) \;,\;  q^n + \lambda_1 q^{n-1/2} + o(q^{n-1/2}) \;\Bigr]
\end{equation}
as $q \to \infty$ through prime powers.
\end{theorem}

\begin{remark}
The interval \eqref{E:interval for large q} is  
a fraction $\lambda_1/(2n)$ of the Hasse--Weil interval, approximately.
\end{remark}

\begin{remark}
For $n=1$, if $q$ is prime, then every integer in $[{q-2q^{1/2}+1} , {q+2q^{1/2}+1}]$
is $\#A(\F_q)$ for some elliptic curve $A$ over $\F_q$.
This fails for $q=p^e$ with $e \ge 2$ because of Remark~\ref{R:Honda-Tate} below.
\end{remark}

\begin{remark}
\label{R:n=2 introduction}
For $n=2$, Theorem~\ref{T:large q} holds if $q$ tends to $\infty$ through primes only.
If instead $q$ tends to $\infty$ through \emph{non-prime} prime powers,
then the constant $\lambda_1 = 2$ (asymptotically 50\% of the Hasse--Weil interval) 
must be replaced by $\lambda_2 \colonequals 4 - 2 \sqrt{2}$ (about 29\% of the Hasse--Weil interval);
see Remark~\ref{R:n=2}.
\end{remark}

\begin{remark}
\label{R:large q ordinary}
If we allow only \emph{ordinary} abelian varieties, then 
Theorem~\ref{T:large q} remains true for $n \ge 3$, as the proof will show,
but for $n=2$ one must use $\lambda_2$ in place of $\lambda_1$, even if $q$ is prime.
\end{remark}

\begin{remark}
DiPippo and Howe proved a result implying that for any $n \ge 2$,
all integers in an interval of the form \eqref{E:interval for large q} with $\lambda_1$ replaced by $1/2$
are realized by ordinary abelian varieties \cite[Theorem~1.4]{DiPippo-Howe-1998}.
Thus Theorem~\ref{T:large q} and Remark~\ref{R:large q ordinary} give an asymptotically optimal improvement of their result.
\end{remark}

Theorem~\ref{T:large q} will be proved in Section~\ref{S:large q}.

\subsection{Effective bounds}
\label{S:introduction effective}

The polynomial constructions we used to prove Theorems~\ref{T:almost-converse} and~\ref{T:large q} are difficult to analyze explicitly for specific values of $q$ and $n$, even when $q=3$.
In Section~\ref{S:effective}, we give \emph{another} construction, and this one, combined with some computations with rigorous error bounds, will allow us to prove the following.

\begin{theorem}
\label{T:effective}
Let $q$ be a prime power.
\begin{enumerate}[\upshape (a)]
\item \label{I:2,3,4,5} 
For each $q \le 5$, every positive integer is $\#A(\F_q)$ for some abelian variety $A$ over $\F_q$.
\item \label{I:constant 3} 
For arbitrary $q$, every integer $\ge q^{3 \sqrt{q} \log q}$ is $\#A(\F_q)$ for some abelian variety $A$ over $\F_q$. 
\end{enumerate}
\end{theorem}

\begin{remark}
Theorem~\ref{T:effective}\eqref{I:2,3,4,5} is best possible: As remarked in \cite{Howe-Kedlaya-preprint}, if $q \ge 7$, then $2$ lies outside the union of the Hasse--Weil intervals \eqref{E:Weil bound}.
\end{remark}

\begin{remark}
Theorem~\ref{T:effective}\eqref{I:constant 3} is best possible too, except for the constant $3$, which we have not attempted to optimize.
It becomes false for large $q$ if $3$ is replaced by any number $\delta < 1/4$,
because if $n = (\delta + o(1)) \sqrt{q} \log q$,
then 
\begin{align*}
    \log \frac{(\sqrt{q}-1)^{2(n+1)}}{(\sqrt{q}+1)^{2n}} 
    &= \log q + o(1) + 2n \log \frac{\sqrt{q}-1}{\sqrt{q} + 1} \\
    &= \log q + o(1) + 2(\delta + o(1)) (q^{1/2} \log q) (- 2 q^{-1/2} + o(q^{-1})) \\
    &= (1 - 4 \delta + o(1)) \log q,
\end{align*}
which means that there is a large gap between the $n$th Hasse--Weil interval and the $(n+1)$st.
\end{remark}

\begin{remark}
\label{R:realizing ordinary}
Suppose that we require $A$ to be ordinary.
Both statements in Theorem~\ref{T:effective} remain true,
except that when $q=4$ one must exclude order~$3$  (that $3$ over $\F_4$ must be excluded
follows from~\cite[Theorem~3.2]{Kadets-2021}).
\end{remark}

\begin{remark}
\label{R:7 and 8}
For $q=7$, the only positive integers not of the form $\#A(\F_7)$ are $2$, $14$, and $17$.
If we require $A$ to be ordinary, then $8$ and $73$ are the only additional integers that must be excluded.
\end{remark}

\begin{remark}
\label{R:realizing squarefree}
Suppose that we require $f_A$ to be squarefree.
Then all the claims in this section remain true except that for $q=7$, the integer $16$ is no longer realized.
\end{remark}

%****************************************************************************
\section{Honda--Tate theory}\label{S:Honda-Tate}

Throughout the paper, if $f$ is a polynomial, then $f^{[i]}$ denotes the coefficient of its degree~$i$ term.
All the results of this section are restatements of results in \cite[Chapter~2]{Waterhouse-1969}.

\begin{theorem}[{Honda--Tate}]
\label{T:Honda-Tate}
A polynomial $f \in \Z[x]$ is the characteristic polynomial
of an ordinary abelian variety $A$ of dimension $n$ over $\F_q$
if and only if
\begin{enumerate}[\upshape \phantom{mm}(a)]
    \item $f$ is monic of degree $2n$; 
    \item \label{I:q-symmetric} $f$ is \defi{$q$-symmetric}, by which we mean $f^{[i]} = q^{n-i} f^{[2n-i]}$ for $i=0,\ldots,n-1$;
    \item \label{I:complex roots} all complex roots of $f$ have absolute value $q^{1/2}$; and
    \item \label{I:middle coefficient} $p \nmid f^{[n]}$.
\end{enumerate}
\end{theorem}

\begin{remark}
Condition~\eqref{I:complex roots} implies \eqref{I:q-symmetric} if $x+q^{1/2}$ and $x-q^{1/2}$ each appear to an even power in the factorization of $f(x)$ over $\C$.
\end{remark}

\begin{remark}\label{R:Honda-Tate}
Let $v \colon \Q_p \to  \Z \union \{\infty\}$ be the $p$-adic valuation.
If in Theorem~\ref{T:Honda-Tate} we replace \eqref{I:middle coefficient} by the weaker condition 
\begin{enumerate}[\upshape \phantom{mm}(d$'$)]
\item \label{I:fancy condition}
the multiplicity $\mu$ of each $\Q_p$-irreducible factor $g$ in $f$ is such that $\mu \, v(g(0))/v(q) \in \Z$,
\end{enumerate}
then we obtain the criterion for $f$ to be the characteristic polynomial of a not-necessarily-ordinary abelian variety $A$ of dimension $n$ over $\F_q$.
If $q$ is prime, then (d$'$) holds automatically.
\end{remark}

%****************************************************************************
\section{Roots on a circle}\label{S:circle}

For $r>0$, let $\C_{\le r}$ be the closed disk $\{z \in \C : |z| \le r\}$.
Let $D = \C_{\le q^{-1/2}}$.  
For 
\begin{align*}
   h(z) &= a_0 + a_1 z + \cdots + a_s z^s \; \in \R[z] \\
\intertext{with $s<2n$, define}
    \widehat{h}(x) &= x^{2n} h(1/x) + q^n h(x/q) \\
    &= a_0 x^{2n} + a_1 x^{2n-1} + \cdots + a_s x^{2n-s} \\
    &\qquad \qquad \qquad {} + q^{n-s} a_s x^s + \cdots + q^{n-1} a_1 x + q^n a_0,
\end{align*}
which is a $q$-symmetric polynomial of degree $\le 2n$
(the notation implicitly depends on a choice of $n$).
To prove Theorem~\ref{T:almost-converse}, 
we will eventually need $\widehat{h}$ for some polynomials $h$ of degree $s>n$, 
in which case the two ranges of exponents of $x$ overlap.

\begin{lemma}\label{L:h to f}
Let $h(z) \in \R[z]$ be a polynomial of degree $<2n$ such that $h$ is nonvanishing on $D$.
Then all complex roots of $\widehat{h}(x)$ have absolute value $q^{1/2}$.
\end{lemma}

\begin{proof}
As $z$ goes around the circle $|z|=q^{-1/2}$, the winding number of $h(z)$ around $0$ is $0$,
so the winding number of $x^n h(1/x)$ as $x$ goes around the circle $|x|=q^{1/2}$ is $n$.
Thus the real-valued function $2 \re(x^n h(1/x)) = x^n h(1/x) + q^n x^{-n} h(x/q)$ on the circle $|x|=q^{1/2}$
crosses $0$ at least $2n$ times.
Multiplying by $x^n$ shows that $\widehat{h}(x)$ has at least $2n$ roots on the circle $|x|=q^{1/2}$.
It cannot have more than $2n$ roots, since $\deg \widehat{h} = 2n$.
\end{proof}

\begin{remark}
If $h(z) = 1 + a_1 z + \cdots + a_n z^n$ with $\sum_{i=1}^n |a_i| q^{-i/2} < 1$,
then $h(D) \subset \{z \in \C: {|z-1|} <1 \}$, so $0 \notin h(D)$.
Thus Lemma~\ref{L:h to f} subsumes \cite[Lemma~3.3.1]{DiPippo-Howe-1998},
which appears also (with a different proof) as \cite[Lemma~2]{Howe-Kedlaya-preprint}.
The feature of Lemma~\ref{L:h to f} that allows us to obtain stronger results 
is that $\{h : 0 \notin h(D)\}$ is closed under multiplication,
a natural property given that one can take products of abelian varieties.
\end{remark}

\begin{remark}
\label{R:squarefree}
The polynomials $\widehat{h}(x)$ produced by Lemma~\ref{L:h to f} are squarefree.
\end{remark}

\begin{remark}
Applying Lemma~\ref{L:h to f} to $h(rx)$ as $r \to 1^-$ shows that the hypothesis
could be weakened to assume only that $h$ is nonvanishing on the \emph{interior} of $D$.
\end{remark}

For use in the proof of Lemma~\ref{L:P^b}, we record the following result.

\begin{lemma}
\label{L:value of R(1)}
Let $R \in \C[z]$ be a polynomial with no zeros inside $D$.
Then
\begin{equation}
   |R(1)| \le q^{(\deg R)/2} |R(1/q)|.
\end{equation}
\end{lemma}

\begin{proof}
By multiplicativity in $R$, we may assume that $R(z) = z-w$ for some $w \in \C$ with $|w| \ge q^{-1/2}$.
We must prove $|(1-w)/(1/q-w)| \le q^{1/2}$.
The M\"obius transformation $M(z) \colonequals (1-z)/(1/q-z)$ maps the circle $|z|=q^{-1/2}$
to a complex-conjugation-invariant circle passing through $M(\pm q^{-1/2}) = \pm q^{1/2}$,
and it maps the exterior to the interior since $M(\infty)=1$.
\end{proof}

%****************************************************************************
\section{Abelian varieties of all sufficiently large orders}\label{S:simple proof}

\begin{theorem}
  \label{T:small range}
  Fix a prime power $q$ and a closed interval $I \subset \R_{>0}$.
  For $n \gg 1$, each integer $m \in q^n I$ is $\#A(\F_q)$
  for some ordinary abelian variety $A$ of dimension $n$ over $\F_q$.
\end{theorem}

\begin{proof}
  For $k \ge 1$, let $\calJ_k$ be the set of power series of the form
  $1+a_k z^k + a_{k+1} z^{k+1} + \cdots$
  with integer coefficients in $[-q/2,q/2]$.
  Choose $k$ such that $1 - \sum_{r \ge k} \lfloor q/2 \rfloor q^{-r/2} \ge 1/2$;
  then $|j(w)| \ge 1/2$ for all $j \in \calJ_k$ and $w \in D$.
  Choose $\epsilon>0$ such that $[1-\epsilon,1+\epsilon] \subset \{j(1/q) : j \in \calJ_k\}$.
  Choose $N$ such that $[(1-\epsilon)^N,(1+\epsilon)^N] \supset I$.
  Then, given $m \in q^n I$, we may choose $j \in \calJ_k$
  with $j(1/q)^N = m/q^n$.
  Write $j^N = h_0 + h_1$
  such that $h_0 \in 1 + z^k \Z[z]$ is of degree $\le n$,
  and $h_1 \in z^{n+1}\Z[[z]]$.
  Let $E = m - \widehat{h}_0(1)$.
  Let 
  \[
     h = h_0 + (E/2) z^n + s (z^{n-1}-((q+1)/2) z^n),
  \]
  where $s \in \{0,1\}$ is chosen so that $p$ does not divide
  the coefficient of $x^n$ in 
  \[
  \widehat{h} = \widehat{h}_0 + E x^n + s (x^{n+1} - (q+1) x^n + q x^{n-1}).
  \]
  Then $\widehat{h}$ is a monic polynomial of degree $2n$ in $\Z[x]$
  and $\widehat{h}(1) = \widehat{h}_0(1) + E = m$.
  The conclusion follows from Lemma~\ref{L:h to f} 
  and Theorem~\ref{T:Honda-Tate} if we can show that $h$ is nonvanishing on $D$.
  We will do so by estimating the error in the approximations
  $h \approx h_0 \approx j^N$.

  Since $j$ has bounded coefficients, induction on $N$ shows that 
  $|(j^N)^{[r]}| = O(r^{N-1})$ as $r \to \infty$,
  uniformly for $j \in \calJ_k$.
  Thus 
  \begin{align*}
     |h_0(1)| &= \biggl\lvert \sum_{r=0}^n (j^N)^{[r]} \biggr\rvert \le \sum_{r=0}^{n} O(r^{N-1}) = O(n^N), \\
     |h_1(1/q)| &= \biggl\lvert \sum_{r=n+1}^\infty (j^N)^{[r]} q^{-r} \biggr\rvert  \le \sum_{r=n+1}^{\infty} O(r^{N-1}) q^{-r} = O(n^{N-1} q^{-n-1}), \\
      |E| &= |m - \widehat{h}_0(1)| = |q^n j(1/q)^N - (q^n h_0(1/q) + h_0(1))|
  	\le |q^n h_1(1/q)| + |h_0(1)| = O(n^N).
\end{align*}
  Now
  \[
  h(z) = j(z)^N - h_1(z) + (E/2) z^n + s (z^{n-1}-((q+1)/2) z^n),
  \]
  so for $w \in D$,
  \[
  |h(w)| \ge 2^{-N} - O(n^{N-1}) q^{-n/2} - O(n^N q^{-n/2}) - O(q \cdot q^{-n/2}) > 0
  \]
  if $n$ is large enough.
\end{proof}

\begin{corollary}
  \label{C:weak existence theorem}
  Fix a prime power $q$.
Every sufficiently large positive integer is $\#A(\F_q)$
for some ordinary abelian variety $A$ over $\F_q$.
\end{corollary}

\begin{proof}
  Apply Theorem~\ref{T:small range} with $I=[1,q]$.
\end{proof}

%****************************************************************************
\section{A congruence condition forcing geometric simplicity and the existence of principal polarizations}
\label{S:congruence condition}

The goal of this section is Proposition~\ref{P:good polynomial mod L}, 
which provides a congruence condition on the characteristic polynomial
of an abelian variety $A$ over $\F_q$ 
which guarantees that $A$ is geometrically simple and isogenous to a principally polarized abelian variety.
Moreover, the congruence condition will be compatible with prescribing $\#A(\F_q)$.
The lemmas in this section are used only to prove Proposition~\ref{P:good polynomial mod L}.

\begin{lemma}
  \label{L:irreducible polynomials}
  For every prime power $q$, prime $\ell \ge 7$ not dividing $q$,
  and integer $n \ge 1$,
  there exists $j(x) \in \F_\ell[x]$ such that
  $j(x)$ and $x^n j(q/x)$ are relatively prime irreducible polynomials
  of degree $n$ not vanishing at $1$.
\end{lemma}

\begin{proof}
  If $n=1$, choose $j(x)=x-a$ where $a \in \F_\ell - \{0,1,q,\pm\sqrt{q}\}$.
  If $n=2$, let $j(x)$ be the minimal polynomial of an element
  $\alpha \in \F_{\ell^2}^\times - \F_\ell^\times$ such that $\alpha \ne q/\alpha$
  and $\alpha^\ell \ne q/\alpha$; there are at least $(\ell^2-\ell) - 2 - (\ell+1) > 0$ such elements $\alpha$.

  Now suppose that $n \ge 3$.
  Let $\alpha$ be a generator of the multiplicative group $\F_{\ell^n}^\times$.
  Let $j(x)$ be the minimal polynomial of $\alpha$ over $\F_\ell$.
  If $j(x)$ and $x^n j(q/x)$ are not relatively prime,
  then $\alpha^{\ell^a} = q/\alpha$ for some $a \in \{0,1,\ldots,n-1\}$.
  Then $\alpha^{(\ell-1)(\ell^a + 1)} = q^{\ell-1} = 1$ in $\F_{\ell^n}$,
  so $\ell^n-1$ divides $(\ell-1)(\ell^a + 1)$,
  contradicting $0 < (\ell-1)(\ell^a + 1) < \ell^n-1$.
\end{proof}

\begin{lemma}
  \label{L:cycle types}
  Let $q$ be a prime power, let $\ell \ge 7$ be a prime not dividing $q$,
  let $n \in \Z_{\ge 1}$, and let $m \in \Z$.
  Suppose that $d_1,\ldots,d_r$ are positive integers summing to $n$
  such that $1$ appears exactly once or twice among $d_1,\ldots,d_r$
  and every other positive integer appears at most once.
  Then there exists a monic $q$-symmetric polynomial $g(x) \in \F_\ell[x]$
  such that 
  \begin{itemize}
      \item $g(1) = m \bmod \ell$,
      \item the roots of $g$ form $n$ \emph{distinct} multiset pairs $\{\alpha,q/\alpha\}$; and
      \item the Frobenius element of $\Gal(\Fbar_\ell/\F_\ell)$ acts on these $n$ pairs 
  as a permutation consisting of cycles of lengths $d_1,\ldots,d_r$.
  \end{itemize}
\end{lemma}

\begin{proof}
  For each $i$ with $d_i \ge 2$, let $j_i(x)$ be the polynomial of degree $d_i$
  provided by Lemma~\ref{L:irreducible polynomials},
  and let $g_i(x) = j_i(x) \cdot x^{d_i} j_i(q/x)$.
  For each $i$ with $d_i=1$, 
  let $g_i(x) = x^2-a_i x+q$ for some $a_i \in \F_\ell$ to be determined.
  Then $g(x) = \prod_{i=1}^r g_i(x)$ gives the correct cycle type,
  and its irreducible factors are distinct,
  except possibly for the factors of the $g_i$ for which $d_i=1$.
  
  If exactly one $d_i$ equals $1$,
  then there is a unique choice of $a_i$ in $\F_\ell$ that makes 
  $g(1) = m \bmod \ell$.
  If $d_i$ and $d_j$ both equal $1$ (with $i \ne j$),
  then there are at least $\ell-1$ choices for $(a_i,a_j)$ that make 
  $g(1) = m \bmod \ell$
  and at most two of these satisfy $a_i=a_j$;
  thus we can ensure $g(1) = m \bmod \ell$
  while making $g$ separable.
\end{proof}

\begin{lemma}
  \label{L:no elliptic curve factor}
  For every prime power $q$, integer $m$, prime $\ell > q+2\sqrt{q}+1$,
  and integer $n \ge 8\sqrt{q}+5$, there exists a monic $q$-symmetric polynomial
  $g(x) \in \F_\ell[x]$ of degree $2n$ such that $g(1)=m \bmod \ell$
  and $g(x)$ has no factor of the form $x^2-\overline{a}x+q$
  with $a \in \Z$ and $|a| \le 2\sqrt{q}$.
\end{lemma}

\begin{proof}
  Since $\ell > q+2\sqrt{q}+1$, none of the polynomials $x^2-\overline{a}x+q$
  vanish at $1\bmod \ell$.
  Lagrange interpolation provides a monic degree~$n$ polynomial
  $j(x) \in \F_\ell[x]$
  such that $j(0)=1$, $j(1)=m$, $j(q)=1$, and $j(\alpha)=1$
  for every root $\alpha \in \Fbar_\ell$ of the quadratic polynomials
  $x^2-\overline{a}x+q$ (the number of values to specify is at most
  $3+2(4\sqrt{q}+1) \le n$).
  Take $g(x) \colonequals j(x) \cdot x^n j(q/x)$.
\end{proof}

\begin{lemma}
  \label{L:generating S_n or S_{n-1}}
  Let $n \ge 3$.
  A subgroup $G$ of $S_n$ containing an $(n-1)$-cycle,
  an $(n-2)$-cycle, and a $2$-cycle
  is either $S_n$ or the stabilizer $S_{n-1}$ of the fixed point
  of the $(n-1)$-cycle.
\end{lemma}

\begin{proof}
 Without loss of generality, the fixed point of the $(n-1)$-cycle is $n$.
  If $G \le S_{n-1}$, then $G$ acts on $\{1,2,\ldots,n-1\}$
  transitively (because of the $(n-1)$-cycle)
  and primitively (because of the $(n-2)$-cycle);
  a primitive subgroup of $S_{n-1}$ containing a $2$-cycle
  is the whole group $S_{n-1}$ \cite[Theorem~8.17]{Isaacs2008}.
  Otherwise $G$ acts on $\{1,\ldots,n\}$ transitively
  (because the orbit of $1$ is larger than $\{1,2,\ldots,n-1\}$)
  and primitively (because of the $(n-1)$-cycle),
  and then the $2$-cycle forces $G = S_n$.
\end{proof}

\begin{lemma}
  \label{L:almost geometrically simple}
  Let $n \ge 5$. 
  Let $A$ be an $n$-dimensional abelian variety over $\F_q$.
  Write $f_A(x) = x^n R(x+q/x)$
  for some monic degree~$n$ polynomial $R(x) \in \Z[x]$.
  If the Galois group of $R$ is $S_n$ or the stabilizer $S_{n-1}$ of a point,
  then $A$ is either geometrically simple
  or a product of geometrically simple abelian varieties over $\F_q$
  of dimensions $n-1$ and $1$.
\end{lemma}

\begin{proof}
  If $A$ is isogenous to $A_1 \times A_2$ over $\F_q$,
  then $R$ factors correspondingly into $R_1 R_2$.
  Since $R$ is either irreducible or a product of irreducible polynomials
  of degrees $1$ and $n-1$,
  the abelian variety $A$ is either simple
  or a product of simple abelian varieties of dimensions $1$ and $n-1$.
  Let $A'$ be the simple factor of dimension $d \in \{n,n-1\}$, 
  and define $R'$ accordingly.

  Suppose that $A'$ is not geometrically simple.
  Let $r>1$ be such that $A'_{\F_{q^r}}$ is not simple.
  Then $f_{A'}$ has roots $\alpha,\beta \in \Qbar$
  giving rise to distinct roots $\alpha+q/\alpha \ne \beta+q/\beta$ of $R'$
  such that $\alpha^r = \beta^r$.
  Now $\beta = \zeta \alpha$ for some root of unity $\zeta$.
  Thus the extension $\Q(\alpha, \zeta) \supset \Q(\alpha+q/\alpha)$, being the compositum of two abelian extensions, is abelian, so its subfield $\Q(\alpha+q/\alpha, \beta+q/\beta)$ is Galois over $\Q(\alpha+q/\alpha)$,
  contradicting the fact that $S_{d-2}$ is not normal in $S_{d-1}$. 
\end{proof}

\begin{lemma}
\label{L:ppav}
For every prime power $q=p^e$, prime $\lambda \ge 7$ such that $q$ is a nonzero square modulo $\lambda$, 
and integers $n \ge 5$ and $m$,
there exists a monic $q$-symmetric degree~$2n$ polynomial $g(x) \in (\Z/\lambda^2 \Z)[x]$
with $g(1) = m \bmod{\lambda^2}$
such that if $A$ is a simple abelian variety over $\F_q$ with $f_A \bmod \lambda^2$ equal to $g$,
then the isogeny class of $A$ contains a principally polarized abelian variety.
\end{lemma}

\begin{proof}
By Hensel's lemma, we can choose $s \in \Z$ such that the discriminant of $x^2-sx+q$ is $0 \bmod \lambda$
but nonzero mod $\lambda^2$.
Replace $s$ by $-s$, if necessary, to make $q+1-s \not\equiv 0 \pmod{\lambda}$.
Choose a monic irreducible polynomial $S(x) \in \F_{\lambda}[x]$ of degree $n-3$.
Choose $a,b \in \F_\lambda$ such that the polynomial $\bar{R} \colonequals (x-s)(x-a)(x-b)S(x) \in \F_\lambda[x]$ is separable and $\bar{R}(q+1) = m \bmod \lambda$; 
this amounts to choosing two elements of $\F_\lambda$ (namely, $q+1-a$ and $q+1-b$) with prescribed product,
not equal to $q+1-s$ or each other, which is possible because $\lambda-1>4$.
Let $R \in (\Z/\lambda^2\Z)[x]$ be a lift of $\bar{R}$ such that $R(s)=0$ and $R(q+1)=m$ in $\Z/\lambda^2\Z$.
Let $g(x) = x^n \, R(x+q/x) \in (\Z/\lambda^2\Z)[x]$.

Suppose that $A$ is a simple abelian variety over $\F_q$ such that $f_A \bmod \lambda^2$ is $g$.
Since $A$ is simple, $f_A$ is a power of an irreducible polynomial \cite[Chapter~2]{Waterhouse-1969},
but its reduction $g \bmod \lambda$ has some simple roots (for example, the roots of $x^{n-3} S(x+q/x)$), 
so $f_A$ must be irreducible, of degree $2n$.
Let $\pi \in \Qbar$ be a root of $f_A$.
Let $K = \Q(\pi)$ and $K^+ = \Q(\pi+q/\pi)$, so $K$ is a CM field and $K^+$ is its totally real subfield.
Since the minimal polynomial of $\pi+q/\pi$ reduces to $\bar{R}$, the extension $K^+/\Q$ is unramified above $\lambda$.
On the other hand, $K/K^+$ is ramified at the prime above $\lambda$ corresponding to the root $s$ of $g$,
because the discriminant of $x^2-sx+q$ has odd valuation $1$.
By \cite[Theorem~1.1]{Howe-1996AG}, the isogeny class of $A$ contains a principally polarized abelian variety.
\end{proof}

\begin{lemma}
\label{L:lambda}
For any prime power $q$, there exists a prime $\lambda$ such that $7 \le \lambda < q^3$ and $q$ is a nonzero square mod $\lambda$.
\end{lemma}

\begin{proof}
We will choose $\lambda$ to be a prime factor of $u^2-q$ for some integer $u$ in $[\sqrt{q}-30,\sqrt{q}+30]$
chosen so that $u^2-q \ne \pm1$ and $u^2-q$ is not divisible by $2$, $3$, or $5$.
There are at least six integers $u$ in $[\sqrt{q}-30,\sqrt{q}+30]$ 
such that $u^2-q$ is not divisible by $2$, $3$, or $5$.
At most two of them satisfy $u^2-q=\pm1$; among the other four are two differing by $30$,
and one of them is prime to $p$.
Thus $u$ can be found.
Then $\lambda \ne 2,3,5,p$, and $\lambda \le (\sqrt{q}+30)^2-q$,
which is less than $q^3$, except for some small $q$ for which we instead compute an explicit $\lambda$.
\end{proof}

\begin{proposition}
  \label{P:good polynomial mod L}
  Given a prime power $q$, there exists a positive integer $L$ 
  such that for any integers $n \gg 1$ and $m$,
  there exists a monic $q$-symmetric polynomial $g(x) \in (\Z/L\Z)[x]$
  of degree $2n$
  with $g(1) = m \bmod L$
  such that any $n$-dimensional abelian variety $A$ over $\F_q$
  whose characteristic polynomial reduces modulo $L$ to $g(x)$
  is ordinary, geometrically simple, and isogenous to a principally polarized abelian variety.
  Moreover, we may choose $L < q^{23}$.
\end{proposition}

\begin{proof}
  Let $\lambda$ be as in Lemma~\ref{L:lambda}.
  Let $L = p \lambda^2 \ell_0 \ell_1 \ell_2 \ell_3$,
  where $p$ is the characteristic, and $p,\lambda,\ell_0,\ldots,\ell_3$ are distinct primes
  such that $\ell_0 > q+2\sqrt{q}+1$ and $\ell_i \ge 7$ for $i=1,\ldots,3$.
  Suppose that $n \ge 8\sqrt{q}+5$.
  Let $\gamma(x) \in \F_p[x]$ be a monic $q$-symmetric polynomial of degree $2n$
  such that $\gamma(1) = m \bmod p$; add $x^{n+1}-x^n$, if necessary, to make $\gamma^{[n]} \ne 0 \bmod{p}$ (here $q$-symmetry means only that $\gamma^{[i]}=0$ for $i<n$).
  Let $g_\lambda(x) \in (\Z/\lambda^2\Z)[x]$ be as in Lemma~\ref{L:ppav}.
  Apply Lemma~\ref{L:no elliptic curve factor} to
  produce a polynomial $g_0(x) \in \F_{\ell_0}[x]$.
  Apply Lemma~\ref{L:cycle types} to produce polynomials
  $g_i(x) \in \F_{\ell_i}[x]$ for $i=1,2,3$ corresponding to the partitions
  \begin{itemize}
  \item $(n-1,1)$
  \item $(n-2,1,1)$
  \item $(n-3,2,1)$ if $n$ is even; and $(n-4,2,1,1)$ $n$ is odd,
  \end{itemize}
  respectively.
  Let $g \in (\Z/L\Z)[x]$ be the monic $q$-symmetric polynomial
  of degree~$2n$ reducing to $\gamma$, the $g_i$, and $g_\lambda$.
  
  Suppose that $A$ is an $n$-dimensional abelian variety over $\F_q$
  such that $f_A(x) \bmod L = g(x)$.
  Write $f_A(x) = x^n R(x+q/x)$.
  Let $G \le S_n$ be the Galois group of $R$,
  which encodes the action of $\Gal(\Qbar/\Q)$
  on the pairs $\{\alpha,q/\alpha\}$ of roots of $F$.
  By choice of $g_1,g_2,g_3$,
  the group $G$ contains permutations $\sigma_1,\sigma_2,\sigma_3$
  whose cycle types are given by the partitions above.
  Raising $\sigma_3$ to an odd power produces a $2$-cycle.
By Lemma~\ref{L:generating S_n or S_{n-1}}, $G$ is $S_n$ or $S_{n-1}$.
By Lemma~\ref{L:almost geometrically simple},
$A$ is either geometrically simple
  or a product of geometrically simple abelian varieties over $\F_q$
  of dimensions $n-1$ and $1$.
  In the second case, $f_A(x)$ would have a factor $x^2-ax+q$
  for some integer $a$ with $|a| \le 2\sqrt{q}$,
  which is ruled out by the choice of $g_0$.
  Thus $A$ is geometrically simple.
  Since $\gamma^{[n]} \ne 0 \bmod p$, $A$ is ordinary.
  By Lemma~\ref{L:ppav}, $A$ is isogenous to a principally polarized abelian variety.
  
  In proving $L < q^{23}$, the worst case is $q=2$,
  in which case we take $L=2 \cdot 7^2 \cdot 11 \cdot 13 \cdot 17 \cdot 19 < 2^{23}$.
\end{proof}

\begin{remark}
\label{R:congruences for p-rank}
It is not hard to adapt Proposition~\ref{P:good polynomial mod L} 
for the purpose of constructing abelian varieties of prescribed order
that have \emph{prescribed $p$-rank}.
Namely, one can prove that it suffices to impose congruences modulo $pq^2$ on the coefficients of 
a $q$-symmetric monic degree~$2n$ polynomial $f$
to guarantee that its Newton polygon is the lowest Newton polygon corresponding to $p$-rank $r$
and that its segments of slope in $[-1/2,0]$ correspond to $\Q_p$-irreducible factors,
in which case the other segments do too by $q$-symmetry,
so that the condition in Remark~\ref{R:Honda-Tate} is satisfied;
moreover one can make these congruences compatible with $f(1) \equiv m \pmod{pq^2}$,
provided that, in the case $r=0$, one has $m \equiv 1 \pmod{p}$.
This last hypothesis is necessary:
if $A$ has $p$-rank $0$, then all roots of $f_A$ have positive $p$-adic valuation, 
so $\#A(\F_q) \equiv 1 \pmod{p}$.
\end{remark}

%****************************************************************************
\section{Chebyshev polynomials}
\label{S:Chebyshev}

Choose the branch of $\sqrt{z^2-1}$ on $\C-[-1,1]$ that is $z + o(1)$ as $z \to \infty$.
Let $\tau(z) = z + \sqrt{z^2-1}$.
Define the $d$th Chebyshev polynomial
\begin{equation}
T_d(z) = \frac{1}{2}\left(\left(z + \sqrt{z^2 - 1}\right)^d + \left(z - \sqrt{z^2 - 1}\right)^d\right) = (\tau(z)^d + \tau(z)^{-d}) / 2.
\end{equation}

\begin{lemma}
\label{L:Chebyshev limit new}
For a suitable choice of $d$th root,
the functions $T_d(z)^{1/d}/z$ and $\tau(z)/z$ extend to holomorphic functions on $\PP^1(\C) \setminus [-1,1]$,
and $T_d(z)^{1/d}/z \to \tau(z)/z$ uniformly on any compact subset of that domain as $d \to \infty$.
\end{lemma}

\begin{proof}
Since $\tau$ is nonvanishing with a simple pole at $\infty$, 
the maximum modulus principle applied to $1/\tau$ 
shows that $|\tau(z)|$ is minimized as $z$ approaches $[-1,1]$, in which case $|\tau(z)| \to 1$.
Thus $|\tau(z)| > 1$ on $\PP^1(\C) \setminus [-1,1]$.
The uniform convergence claim now follows from $T_d(z)/z^d = \tfrac12z^{-d}(\tau(z)^d + \tau(z)^{-d})$.
\end{proof}

\begin{proposition}
\label{P:construction of P}
Let $I$ be a closed interval contained in $\Ia$ $($see~\eqref{E:I attainable}$)$.
Then for even $d \gg 1$, there exists a degree~$d$ polynomial $P(z) \in \R[z]$ such that
\begin{enumerate}[\upshape \phantom{mm}(a)]
    \item $P(0)=1$;
    \item \label{I:positive} $P$ is positive on $\R$; 
    \item \label{I:lower bound on D} $|P(w)|^{1/d} \ge q^{-1/4}$ for all $w \in D \colonequals \C_{\le q^{-1/2}}$; and
    \item \label{I:range for 1/q value} $(q P(1/q)^{2/d}, q P(-1/q)^{2/d})$ contains $I$.
\end{enumerate}
\end{proposition}

\begin{remark}
In Appendix~\ref{A:potential theory}, we use potential theory to prove that Proposition~\ref{P:construction of P} is optimal in the sense that it fails if $\Ia$ is enlarged.
\end{remark}

\begin{proof}
For $\epsilon>0$ to be specified later, let 
\begin{align*}
    \ell(z) &= (q^{1/2}/2) z - (q^{1/2}-1), \\
    f_d(z) &= 2 q^{-d/4} z^{d/2} \, T_{d/2}(\ell(z+1/z)), \\ 
    P(z) &= f_d((1-\epsilon)q^{1/2} z).
\end{align*}
\begin{enumerate}[\upshape (a)]
\item
The leading coefficient of $T_{d/2}$ is $2^{d/2-1}$,
so $f_d(0) = 2 q^{-d/4} 2^{d/2-1} (q^{1/2}/2)^{d/2} = 1$ and $P(0)=f_d(0)=1$.

\item
The roots of $T_{d/2}$ are in $[-1,1)$,
and $\ell^{-1}([-1,1)) \subset (-2,2)$,
so all the roots of $f_d(z)$ are on the unit circle and not at $\pm 1$.
Thus $f_d$ does not change sign on $\R$.
Since $f_d(0)>0$, the sign is positive.
Thus $P$ is positive on $\R$.

\item
The function $(1-\epsilon)q^{1/2}z$ maps $D$ to $\C_{\le 1-\epsilon}$,
so we need to prove that $|f_d|^{1/d} \ge q^{-1/4}$ on $\C_{\le 1-\epsilon}$.
First, $z T_{d/2}(\ell(z+1/z))^{2/d}$ is the product of
the polynomial $z \, \ell(z+1/z)$ and holomorphic function $T_{d/2}(\ell(z+1/z))^{2/d}/\ell(z+1/z)$
on $\C_{\le 1-\epsilon}$,
so Lemma~\ref{L:Chebyshev limit new} implies that 
\begin{equation}
\label{E:f_d}
   |f_d(z)|^{1/d} \to q^{-1/4} |z|^{1/2} \left|\tau(\ell(z+1/z)) \right|^{1/2}
\end{equation}
uniformly for $z \in \C_{\le 1-\epsilon}$.
The function $z \, \tau(\ell(z+1/z))$ is holomorphic, nonconstant, and nonvanishing on $\C_{<1}$,
and it extends to a continuous function on $\C_{\le 1}$ having absolute value $\ge 1$ on the boundary,
so the maximum modulus principle applied to its inverse
shows that there exists $M>1$ such that $|z \, \tau(\ell(z+1/z))| > M$ for all $z \in \C_{\le 1-\epsilon}$.
The lower bound on $|f_d|$ follows for $d \gg 1$.

\item 
It suffices to prove that $\lim_{\epsilon \to 0^+} \lim_{d \to \infty} q P(1/q)^{2/d}$ equals
the left endpoint of $\Ia$, and likewise at the other end.
In fact, \eqref{E:f_d} implies that $\lim_{d \to \infty} qP(1/q)^{2/d}$
is a continuous function of $\epsilon \in [0,1]$, so we may simply \emph{substitute} $\epsilon=0$.
Then 
\begin{align*}
    \lim_{d \to \infty} q P(1/q)^{2/d} &= \lim_{d \to \infty} q f_d(q^{-1/2})^{2/d} \\
    &= q \cdot q^{-1/2} q^{-1/2} \,  |\tau(\ell(q^{-1/2} + q^{1/2}))| \\
    & = \tau(q/2 - q^{1/2} + 3/2).
\end{align*}
Similarly, $\lim_{\epsilon \to 0^+} \lim_{d \to \infty} q P(-1/q)^{2/d} = |\tau(-q/2 - q^{1/2} + 1/2)| = \tau(q/2 + q^{1/2} - 1/2)$.\qedhere
\end{enumerate}
\end{proof}

%****************************************************************************
\section{Construction of polynomials}
\label{sec:range}

We now begin the proof of Theorem~\ref{T:almost-converse}.
Let $I$ be a closed interval in $\Ia$.
Let $P(z)$ be as in Proposition~\ref{P:construction of P} and let $d = \deg P$;
we may assume that $d \ge 53$.

The polynomial $P$ was optimized to have a small value at $1/q$ and large value at $-1/q$.
Lemma~\ref{L:P^b} below shows that this makes $\reallywidehat{P^b}(1)$ small and $\reallywidehat{P(-z)^b}(1)$ large,
where $b$ is chosen to make $P^b$ of degree close to $2n$.
The polynomial $Q$ in Lemma~\ref{L:Q} interpolates between $P(z)$ and $P(-z)$ 
to make $\reallywidehat{Q^b}(1)$ equal a prescribed intermediate value.

\begin{lemma}
\label{L:P^b}
Let $b=b(n)$ and $\ell=\ell(n)$ be functions of $n$ tending to $\infty$
such that $\deg P^b = 2n-2\ell$ and $\ell = o(n)$.
Then 
\begin{equation*}
     \reallywidehat{P^b}(1)^{1/n} \To q P(1/q)^{2/d} \qquad \textup{and} \qquad
     \reallywidehat{P(-z)^b}(1)^{1/n} \To  q P(-1/q)^{2/d}
\end{equation*}
 as $n \to \infty$.  $($Recall that $\reallywidehat{P^b}(1) \colonequals q^nP^b(1/q) + P^b(1)$, which depends on $n$.$)$
\end{lemma}

\begin{proof}
We have 
\[
   \reallywidehat{P^b}(1) = q^n P^b(1/q) + P^b(1) = (q^n + O(q^{n-\ell})) P^b(1/q)
\]
by Lemma~\ref{L:value of R(1)} applied to $P^b$.
Taking $n$th roots yields the left endpoint limit, since $\ell \to \infty$ and $b/n = (2n-2\ell)/(dn) \to 2/d$. 
The right endpoint limit follows similarly.
\end{proof}

Choose integers $\ell = \ell(n)$ and $b = b(n)$ such that $\ell = 4 \log_q n + O(1)$ and $bd = 2n - 2\ell$.
The statements in the rest of this section will hold if $n$ is sufficiently large.
Given $m \in \Z$ such that $m^{1/n} \in I$, we want to construct an $n$-dimensional, ordinary, geometrically simple, principally polarized abelian variety $A$ with $\#A(\F_q)=m$.

\begin{lemma}
\label{L:Q}
There exists $Q(z) \in 1+ z \R[z]$ of degree $\le d$ such that $Q$ is positive on $\R$,
$\reallywidehat{Q^b}(1) = m$, 
and $|Q(w)|^{1/d} \ge q^{-1/4}$ for all $w \in D$.
\end{lemma}

\begin{proof}
Because $n$ is sufficiently large,
Proposition~\ref{P:construction of P}\eqref{I:range for 1/q value} and Lemma~\ref{L:P^b} 
show that
\begin{equation}
\label{I:interval containment}
   (\reallywidehat{P^b}(1)^{1/n} , \reallywidehat{P(-z)^b}(1)^{1/n}) \;\supset \; I \; \ni m^{1/n}.
\end{equation}
By the intermediate value theorem, 
there exists $s \in [-1,1]$  
such that the polynomial 
\[
   Q(z) \colonequals P(sz) \in 1 + z \R[z]
\]
satisfies $\reallywidehat{Q^b}(1)^{1/n} = m^{1/n}$.
Thus $\reallywidehat{Q^b}(1) = m$.
Moreover, $Q$ is positive on $\R$, and $|Q(w)|^{1/d} = |P(sw)|^{1/d} \ge q^{-1/4}$ for all $w \in D$ by Proposition~\ref{P:construction of P}(\ref{I:positive},\ref{I:lower bound on D}).
\end{proof}

In the rest of this section, the implied constant in big-$O$ notation may depend on $q$, $L$, $d$, $P$, and $Q$,
but not on $n$.

The polynomial $Q^b$ has real coefficients.
We could round them to the nearest integer to produce a polynomial $h \in \Z[x]$
and adjust the middle coefficients to make $\widehat{h}(1) = m$, as in Section~\ref{S:simple proof}, 
but it turns out that we cannot guarantee that such an $h$ is nonvanishing on $D$,
as required for Lemma~\ref{L:h to f}.
So instead we adjust the coefficients of $Q$ (inside the $b$th power) by only $O(1/n)$ each
to make the first $d$ coefficients of $\reallywidehat{Q^b}$ integral
(and to make them satisfy the congruences in Proposition~\ref{P:good polynomial mod L}),
and then, to correct the later coefficients, we add correction polynomials designed to be small on $D$,
because as we go along, we need to bound the difference between $Q^b$ and the final $h$ 
to ensure that $h$ is still nonvanishing on $D$.

Let us outline the entire construction; then in a series of lemmas, we will prove that the steps make sense.
\begin{construction}
\label{C:best construction}
\hfill
\begin{enumerate}[\upshape 1.]
    \item Let $Q \in 1 + z \R[z]$ be as in Lemma~\ref{L:Q}.
    \item Let $g \in (\Z/L\Z)[x]$ be as in Proposition~\ref{P:good polynomial mod L}.
    \item Let $Q_0 = Q$.
    \item\label{I:a_i} For $i=1,\ldots,d-1$ in turn, let $a_i \in [0,L/b)$ and $Q_i \colonequals Q_{i-1} + a_i z^i$ and $h_i \colonequals Q_i^b$ be such that $\widehat{h}_i^{[2n - i]} \in \Z$ and $\widehat{h}_i^{[2n - i]} \equiv g^{[2n-i]} \pmod{L}$.
    \item \label{I:c} Let $\Qtilde = Q_{d-1} - c z^d$ and $h_d = \Qtilde^b$, where $c \in \R$ is chosen so that $\widehat{h}_d(1)=m$.
    \item Define ``correction polynomials'' as follows:
    \begin{itemize}
    \item For $i=d,\ldots,\ell-1$, let $\kk_i = z^i \, \Qtilde(z)^b$.  
    \item For $i=\ell,\ldots,n-1$, let $\kk_i= z^i \, \Qtilde(z)^a$, where $a \in \Z_{\ge 0}$ is chosen as large as possible such that $\deg \kk_i < 2n-i$.
    \item Define $\kk_n = z^n/2$.
    \end{itemize}
    The definitions are so that $\widehat{\kk}_i$ is monic of degree $2n-i$ for all integers $i \in [d,n]$.
    \item \label{I:r_i} For $i=d,\ldots,n-1$, let $r_i \in [0,L)$ and $s_{i+1} \in \R_{\ge 0}$ and $h_{i+1} \colonequals h_i + r_i \kk_i - s_{i+1} \kk_{i+1}$, where $r_i$ is such that $\widehat{h}_{i+1}^{[2n-i]} \in \Z$ and $\widehat{h}_{i+1}^{[2n-i]} \equiv g^{[2n-i]} \pmod{L}$, and $s_{i+1}$ is such that $\widehat{h}_{i+1}(1)=m$.
    \item Let $A$ be an abelian variety over $\F_q$ with $f_A=\widehat{h}_n$.
\end{enumerate}
\end{construction}

\begin{lemma}\label{L:a_i}
The $a_i$ can be chosen as specified in Step~\ref{I:a_i}, and they are $O(1/n)$.
\end{lemma}

\begin{proof}
In Step~\ref{I:a_i},
once $a_1,\ldots,a_{i-1}$ have been fixed, 
$\widehat{h}_i^{[2n-i]}$ as a function of $a_i$ is a linear polynomial with leading coefficient $b$,
so $a_i \in [0,L/b)$ can be found.
Then $a_i = O(L/b) = O(1/n)$.
\end{proof}

\begin{lemma}
\label{L:c}
\hfill
\begin{enumerate}[\upshape (a)]
\item
The real number $c$ can be chosen as specified in Step~\ref{I:c}, and $c$ is $O(1/n)$.
\item \label{I:Qtilde positive} We have $\Qtilde(1)>0$ and $\Qtilde(1/q)>0$.
\item  \label{I:Qtilde bounded} The values $\Qtilde(1)$ and $\Qtilde(1/q)$ are $O(1)$.
\end{enumerate}
\end{lemma}

\begin{proof}\hfill
\begin{enumerate}[\upshape (a)]
\item
Since $a_i \ge 0$, we have $Q_{d-1} \ge \cdots \ge Q_0 = Q > 0$ on $\R_{\ge 0}$,
so 
\begin{equation}\label{E:Q_{d-1}^b}
    \reallywidehat{Q_{d-1}^b}(1) \ge \reallywidehat{Q^b}(1) = m.
\end{equation}
Let
\[
c' \colonequals q^{d-1} a_1 + q^{d-2} a_2 + \ldots + q a_{d-1}.
\]
Let $R = Q_{d-1} - c' z^d$.
Then 
\[
R(1) = Q_{d-1}(1) - c' = Q(1) - (q^{d-1} - 1) a_1 - \cdots - (1-1) a_{d-1} \; \in \; (0,Q(1)],
\]
for large $n$, by Lemma~\ref{L:a_i}, and
\begin{align*}
   R(1/q) &= Q_{d-1}(1/q) - c'/q^d \\
   &= (Q(1/q) + a_1 q^{-1} + \cdots +a_{d-1} q^{-(d-1)}) - (a_1 q^{-1} + \cdots +a_{d-1} q^{-(d-1)}) \\
   &= Q(1/q) > 0,
\end{align*}
so
\begin{equation}\label{E:R^b}
   \reallywidehat{R^b}(1)  \le \reallywidehat{Q^b}(1) = m.
\end{equation}
By \eqref{E:Q_{d-1}^b} and \eqref{E:R^b} and the intermediate value theorem,
there exists $c \in [0,c']$ such that $\reallywidehat{(Q_{d-1}-c z^d)^b}(1)= m$.
Moreover, $c=O(c') = O((d-1) q^{d-1} (1/n)) = O(1/n)$.
\item We have $\Qtilde(1) \ge R(1) > 0$ and $\Qtilde(1/q) \ge R(1/q) > 0$.
\item For $w \in \{1,1/q\}$, we have $\Qtilde(w) = Q(w) + O(1/n)$, and $Q(w) \in P([-1,1])$,
an interval independent of $n$.\qedhere
\end{enumerate}
\end{proof}

Lemmas \ref{L:Qtilde lower bound} through~\ref{L:h_n on disk} show that $\Qtilde^b$ is large enough on $D$
and the corrections are small enough that $h_n$ is nonvanishing on $D$.

\begin{lemma}
\label{L:Qtilde lower bound}
We have $|\Qtilde(w)| \ge q^{-d/4} - O(1/n)$ for every $w \in D$.
\end{lemma}

\begin{proof}
By Lemma~\ref{L:Q}, $|Q(w)| \ge q^{-d/4}$,
and $|\Qtilde(w)|$ differs from $|Q(w)|$ by at most $|a_1 w + \cdots + a_{d-1} w^{d-1} - c w^d| = O(1/n)$,
by Lemmas \ref{L:a_i} and~\ref{L:c}.
\end{proof}

\begin{lemma}
\label{L:k_i}
We have $\kk_i(1)>0$ and $\kk_i(1/q)>0$. 
\end{lemma}

\begin{proof}
These follow from Lemma~\ref{L:c}\eqref{I:Qtilde positive}. 
\end{proof}

\begin{lemma}
\label{L:r_i}
The $r_i \in [0,L)$ and $s_{i+1}$ can be chosen as specified in Step~\ref{I:r_i}, 
and $s_{i+1}$ is $O(1)$.
For $i=d,\ldots,\ell-2$, we have the more precise bound $s_{i+1} \in [0,qL]$.
\end{lemma}

\begin{proof}
This is similar to the proof of Lemma~\ref{L:c}.
The value $\widehat{h}_{i+1}^{[2n-i]}$ is $r_i$ plus terms that have already been fixed,
so there is a unique choice $r_i \in [0,L)$ such that 
$\widehat{h}_{i+1}^{[2n-i]} \in \Z$ and $\widehat{h}_{i+1}^{[2n-i]} \equiv g^{[2n-i]} \pmod{L}$.

We seek $s_{i+1}$ making the value
$\widehat{h}_{i+1}(1) = m + r_i \, \widehat{\kk}_i(1) - s_{i+1} \, \widehat{\kk}_{i+1}(1)$
equal to $m$.
By Lemma~\ref{L:k_i}, 
\begin{equation}
\label{E:ge}
    m + r_i \, \widehat{\kk}_i(1) \ge m.
\end{equation}
Let $V=\kk_i/\kk_{i+1}$ and $v = \max\{V(1),V(1/q)\}$.
By Lemma~\ref{L:k_i}, $\widehat{\kk}_i(1) \le v \, \widehat{\kk}_{i+1}(1)$, so
\begin{equation}
\label{E:le}
    m + r_i \, \widehat{\kk}_i(1) - v r_i \, \widehat{\kk}_{i+1}(1) \le m.
\end{equation}
Now \eqref{E:ge}, \eqref{E:le}, and the intermediate value theorem
yield $s_{i+1} \in [0,v r_i] \subseteq [0, v L]$
making $\widehat{h}_{i+1}(1)=m$.

To bound $s_{i+1}$, we need to bound $v$.
The function $V$ is $1/z$, $\Qtilde(z)/z$, or $2/z$;
accordingly, $v$ is $q$, $O(1)$, or $2q$,
with the middle case following from Lemma~\ref{L:c}(\ref{I:Qtilde positive},\ref{I:Qtilde bounded}).
In every case, $v = O(1)$, so $s_{i+1} = O(1)$.
If $i \in [d,\ell-1)$, then $V=1/z$, so $v=q$, so $s_{i+1} \in [0,qL]$.
\end{proof}

\begin{lemma}
\label{L:h_n on disk}
The polynomial $h_n$ is nonvanishing on $D$.
\end{lemma}

\begin{proof}
By construction,
\[
   h_n = \Qtilde^b + \sum_{i=d}^{n-1} (r_i \kk_i - s_{i+1} \kk_{i+1}),
\]
so it suffices to prove that 
\begin{equation}
\label{E:less than 1}
    \sum_{i=d}^{n-1} \left| \frac{r_i \kk_i}{\Qtilde^b} \right| + \sum_{i=d+1}^{n} \left| \frac{s_i \kk_i}{\Qtilde^b} \right| < 1
\end{equation}
on $D$.
We claim that 
\begin{equation}
\label{E:k_i/Qtilde^b}
   \left| \frac{\kk_i}{\Qtilde^b} \right| \le 
   \begin{cases}
       q^{-i/2} & \textup{ if $i \in [d,\ell)$,} \\
       O(n^{-2}) & \textup{ if $i \in [\ell,n]$,}
   \end{cases}
\end{equation}
on $D$.
The case $i \in [d,\ell)$ follows since $\kk_i/\Qtilde^b=z^i$.
In particular, for $i \in [\ell-d/2,\ell)$, we have
$|\kk_i/\Qtilde^b| \le q^{-(\ell-d/2)/2} = O(q^{-\ell/2}) = O(n^{-2})$.
From then on, changing $i$ to $i+d/2$ multiplies $|\kk_i/\Qtilde^b|$
by $|z^{d/2}/\Qtilde| \le q^{-d/4} / (q^{-d/4} - O(1/n)) = 1 + O(1/n)$ 
by Lemma~\ref{L:Qtilde lower bound} (or, at the last step with $i+d/2=n$, by $|(z^n/2)/z^i| = |z^{d/2}/2| \le 1$);
this happens fewer than $n$ times, and $(1+O(1/n))^n =O(1)$, 
so \eqref{E:k_i/Qtilde^b} for $i \in [\ell,n]$ follows.

By Lemma~\ref{L:r_i} and \eqref{E:k_i/Qtilde^b}, the left hand side of \eqref{E:less than 1} is at most 
\[
    \sum_{i=d}^{\ell-1} L q^{-i/2} + \sum_{i=\ell}^{n-1} L \, O(n^{-2}) + \sum_{i=d+1}^{\ell-1} qL \, q^{-i/2} + \sum_{i=\ell}^n O(1) \, O(n^{-2}) \le \frac{2 L q^{-(d-1)/2}}{1-q^{-1/2}} + O(1/n) < 1
\]
if $n$ is large, since $L < q^{23}$ and $d \ge 53$.
\end{proof}

\begin{lemma}
\label{L:hhat_n}
The polynomial $\widehat{h}_n$ is monic of degree $2n$.
Also, $\widehat{h}_n \in \Z[x]$ and $\widehat{h}_n \equiv g \pmod{L}$.
\end{lemma}

\begin{proof}
In Steps \ref{I:a_i} and~\ref{I:r_i}, adjusting $h_i$ to produce $h_{i+1}$ does not change 
the coefficients of $z^{2n}$, $z^{2n-1}$, \dots, $z^{2n-i}$ in $\widehat{h}_i$,
which are integers congruent modulo $L$ to the corresponding coefficients of $g$;
by $q$-symmetry, the same holds for the coefficients of $1$, $z$, \dots, $z^i$.
Thus $\widehat{h}_n$ is monic and has integer coefficients congruent to the coefficients of $g$,
except perhaps the coefficient of $z^n$;
actually it holds for this coefficient too
since $\widehat{h}_n(1)$ is an integer (namely, $m$) and $\widehat{h}_n(1) = m \equiv g(1) \pmod{L}$.
\end{proof}

\begin{proof}[End of proof of Theorem~\ref{T:almost-converse}]\hfill
\begin{itemize}
\item The polynomial $\widehat{h}_n$ is monic of degree $2n$, with integer coefficients, by Lemma~\ref{L:hhat_n}.
\item The polynomial $\widehat{h}_n$ is $q$-symmetric, by definition of the hat.
\item All complex roots of $\widehat{h}_n$ have absolute value $q^{1/2}$, by Lemmas \ref{L:h_n on disk} and~\ref{L:h to f}.
\item The characteristic $p$ does not divide $\widehat{h}_n^{[n]}$, because by Lemma~\ref{L:hhat_n}, $\widehat{h}_n^{[n]}$ is congruent modulo $L$ to $g^{[n]}$, which is nonzero modulo $p$, and $p \mid L$, by construction of $g$.
\end{itemize}
By Theorem~\ref{T:Honda-Tate}, there exists an ordinary $n$-dimensional abelian variety $A$ over $\F_q$
with $f_A = \widehat{h}_n$.
Then $\#A(\F_q)=f_A(1) = \widehat{h}_n(1) = m$.
By Proposition~\ref{P:good polynomial mod L}, $A$ is geometrically simple,
and principally polarized after replacing $A$ by an isogenous abelian variety.
\end{proof}

%****************************************************************************
\section{Large \texorpdfstring{$q$}{q} limit} 
\label{S:large q}

\begin{figure}[H]
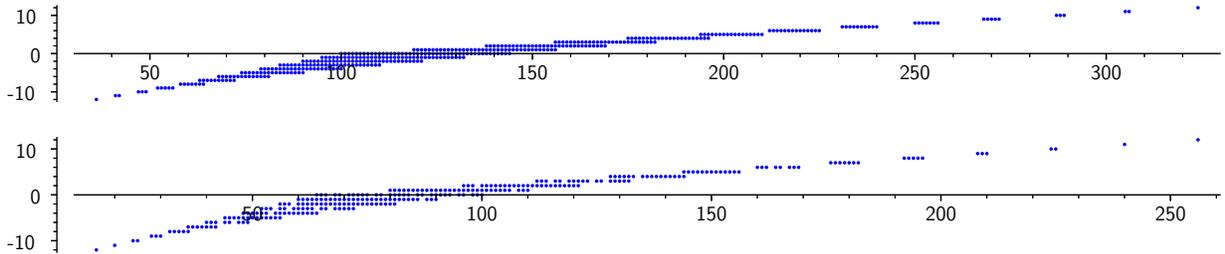

\centering
\begin{subfigure}{\textwidth}
\resizebox{\textwidth}{!}{\input{q11.pgf}}
%\fbox{\parbox{\textwidth}{This image is heavy, check q11.pdf}}
\end{subfigure}

\begin{subfigure}{\textwidth}
\resizebox{\textwidth}{!}{\input{q9.pgf}}
%\fbox{\parbox{\textwidth}{This image is heavy, check q9.pdf}}
\end{subfigure}
\caption{For $q=11$ and $q=9$, respectively, the graph shows all the points $(\#A(\F_q),c_1)$, where $A$ ranges over abelian surfaces over $\F_q$, and $c_1 = G^{[n-1]} = f_A^{[2n-1]}$ with $n=2$; see~\eqref{E:G coefficients}.}
\label{F:11 and 9}
\end{figure}

In this section, we prove Theorem~\ref{T:large q},
which for fixed $n$ and large $q$ determines the largest subinterval of the Hasse--Weil interval 
in which all integers are realizable as $\#A(\F_q)$ for an $n$-dimensional abelian variety $A$ over $\F_q$.
First let us explain the idea.
For any $n$-dimensional abelian variety $A$ over $\F_q$, we have $f_A(x) = x^n \, G(x+q/x)$
for some polynomial 
\begin{equation}
\label{E:G coefficients}
    G(x) = x^n + c_1 x^{n-1} + c_2 x^{n-2} + \cdots + c_n \in \Z[x]
\end{equation}
all of whose roots lie in $[-2q^{1/2},2q^{1/2}]$.
Then $c_i = O(q^{i/2})$, and 
\[
   \#A(\F_q) = f_A(1) = G(q+1) = (q+1)^n + c_1 (q+1)^{n-1} + c_2 (q+1)^{n-2} + \cdots + c_n.
\]
For each integer $c_1$ in the possible range $[-2n q^{1/2},2n q^{1/2}]$,
let $I_{c_1}$ be the smallest interval containing the possible values of
$c_2 (q+1)^{n-2} + \cdots + c_n$; then we prove that the ranges for $c_2$, \dots, $c_n$ are large enough that
all integers in $I_{c_1}$ are realized, possibly ignoring a negligible fraction of the interval at the ends.
The interval $I_{c_1}$ has width $O(q^{n-1})$ and does not change much when $c_1$ is incremented by $1$ --- its endpoints
move by $o(q^{n-1})$.
The big-$O$ constant matters: for $c_1$ close to the extremes of its range
(with $|c_1|$ greater than about $\left(2n -\sqrt{\frac{2n}{n-1}} \right) q^{1/2}$),
it turns out that $I_{c_1}$ has length significantly less than $q^{n-1}$,
so that there is a gap between the intervals $(q+1)^n + c_1(q+1)^{n-1} + I_{c_1}$
and $(q+1)^n + (c_1+1)(q+1)^{n-1} + I_{c_1+1}$, a gap in which $\#A(\F_q)$ cannot lie; see Lemma~\ref{L:unrealizable}.
On the other hand, for the $c_1$ towards the middle of the range,
$I_{c_1}$ has width significantly greater than $q^{n-1}$, 
so the intervals $(q+1)^n + c_1(q+1)^{n-1} + I_{c_1}$ overlap to cover a large interval
in the middle of the Hasse--Weil interval.
Figure~\ref{F:11 and 9} shows these overlapping intervals when $n=2$ and $q \in \{11,9\}$;
for the non-prime $9$, there is an additional phenomenon explained in Remark~\ref{R:n=2}.

As the previous paragraph indicates, the coefficients of $x^{n-1}$ and $x^{n-2}$ are what matter most.
After using the normalization $g(x) \colonequals q^{-n/2} G(q^{1/2} x)$,
we are led to study
\begin{align*}
  \calG &\colonequals \{ \, g \in \R[x] \colon \textup{ $g$ is monic of degree $n$ with all roots in $[-2,2]$} \, \} \\
  \calS &\colonequals \{ \, (g^{[n-1]},g^{[n-2]}) \in \R^2 : g \in \calG \, \}.
\end{align*}
Let $\lambda_1 = 2n - \sqrt{\tfrac{2n}{n-1}}$ and $\lambda_2 = 2n - \sqrt{\tfrac{4n}{n-1}}$.

\begin{lemma}
\label{L:B_diff}
If $n \ge 2$, then there exist continuous functions $B_{\min},B_{\max} \colon [-2n,2n] \to \R$ such that 
\begin{enumerate}[\upshape (a)]
\item \label{I:B_min,B_max} We have $\calS = \{\, (a,b) \in [-2n,2n] \times \R : B_{\min}(a) \le b \le B_{\max}(a) \, \}$.
\item \label{I:B_diff} The difference $B_{\diff}(a) \colonequals B_{\max}(a) - B_{\min}(a)$ is 
\begin{itemize}
\item nonnegative on $[-2n,2n]$, positive on $(-2n,2n)$, 
\item less than $1$ if $\lambda_1 < |a| \le 2n$, greater than $1$ if $|a|<\lambda_1$, 
\item less than $2$ if $\lambda_2 < |a| \le 2n$, and greater than $2$ if $|a|<\lambda_2$.
\end{itemize}
\item \label{I:g_{a,b}} There exists a compact subset $\calG_0 \subset \calG$ surjecting onto $\calS$
such that any $g \in \calG_0$ mapping into the interior of $\calS$ has distinct roots in $(-2,2)$.
\end{enumerate}
\end{lemma}

\begin{proof}
If $g=\prod_{i=1}^n (x-r_i)$, then $(g^{[n-1]},g^{[n-2]}) = (-\sum r_i,\sum_{i<j} r_i r_j)$.
Given $a \in [-2n,2n]$, let 
\[
    \calC_a = \{ (r_1,\ldots,r_n) \in [-2,2]^n : \textstyle \sum r_i =-a\}.
\]
Since $\calC_a$ is compact and connected, \eqref{I:B_min,B_max} holds
with $B_{\min}$ and $B_{\max}$ being the minimum and maximum of $\sum_{i<j} r_i r_j$
on $\calC_a$.
If any two of the $r_i$ are different, then $\sum_{i<j} r_i r_j$ can be increased
by bringing them closer together;
thus the maximum occurs when the $r_i$ are all equal,
so $B_{\max}(a) = \binom{n}{2} (a/n)^2$.
If there are two $r_i$ in $(-2,2)$, then $\sum_{i<j} r_i r_j$ can be decreased 
by moving them slightly apart;
thus the minimum occurs when all but one $r_i$ are at $\pm 2$.
Given $a$, there is at most one such $(r_1,\ldots,r_n)$ with $\sum r_i=-a$ up to permuting coordinates --- 
as $a$ increases, the roots move linearly from $2$ to $-2$ one at a time,
so $B_{\min}$ is the piecewise-linear continuous function such that for each $k \in \{0,\ldots,n-1\}$,
\[
   B_{\min}(a) = (4k-2n+2) a -8 k^2+8 k (n-1)-2 (n-1) n \quad\textup{for $a \in [4k-2n,4k-2n+4]$}.
\]
The minimum value of $B_{\diff}$ on $[4k-2n,4k-2n+4]$ is 
\[
   B_{\diff}(4k-2n + 4k/(n-1)) = 8k(n-1-k)/(n-1),
\]
which for $k \in \{1,\ldots,n-2\}$ is at least $8(n-2)/(n-1) \ge 4$.
On the other hand, for $t \in [0,4]$, we have $B_{\diff}(2n-t) = B_{\diff}(-2n+t) = \frac{n-1}{2n} t^2$.
The claims in \eqref{I:B_diff} follow.

Given $a$, let $\prod_{i=1}^n (x-r_i)$ and $\prod_{i=1}^n (x-r_i'')$
be the polynomials realizing $B_{\min}(a)$ and $B_{\max}(a)$, 
each with roots listed in increasing order.  (So all but one $r_i$ are $\pm2$, and $r_i''=-a/n$ for all $i$.)
Let $\epsilon \ge 0$ be the distance from $-a/n$ to the boundary of $[-2,2]$,
and let $r_1',\ldots,r_n'$ be an arithmetic progression with $r_1' = -a/n-\epsilon/2$ and $r_n'=-a/n+\epsilon/2$.
For each $s \in [0,1]$, consider the monic degree~$n$ polynomial whose roots are $(1-s) r_i + s r_i'$ for $i=1,\ldots,n$
and the analogous polynomial with roots $(1-s) r_i' + s r_i''$.
These depend continuously on $(a,s) \in [-2n,2n] \times [0,1]$,
so the set of all such polynomials is a compact subset $\calG_0$ of $\calG$.
For fixed $a$, the coefficients of $x^{n-2}$ in these polynomials vary continuously from $B_{\min}(a)$ to $B_{\max}(a)$,
so $\calG_0 \to \calS$ is surjective.
Finally, by construction, all polynomials in $\calG_0$ except for the ones realizing $B_{\min}(a)$ and $B_{\max}(a)$
have distinct roots in $(-2,2)$.
\end{proof}

\begin{lemma}
\label{L:unrealizable}
Suppose $n \ge 2$.
For $\lambda \in \R$ satisfying $\lambda_1 < |\lambda| < 2n$,
there exists $\epsilon>0$
such that if $q$ is sufficiently large
and $r = \lfloor \lambda q^{1/2} \rfloor$,
then the interval
\begin{equation}
\label{E:unrealizable for large q}
\bigl[ \; (q+1)^n + r (q+1)^{n-1} + (B_{\max}(\lambda) + \epsilon) q^{n-1} \; , \; (q+1)^n + (r+1) (q+1)^{n-1} + (B_{\min}(\lambda) - \epsilon) q^{n-1} \; \bigr]
\end{equation}
has width $>1$ and does not contain $\#A(\F_q)$ for any $n$-dimensional abelian variety $A$ over $\F_q$.
\end{lemma}

\begin{proof}
By Lemma~\ref{L:B_diff}\eqref{I:B_diff}, $B_{\diff}(\lambda) < 1$.
Choose $\epsilon>0$ such that $B_{\diff}(\lambda) < 1-2\epsilon$.
Then the width of the interval \eqref{E:unrealizable for large q} is 
$(q+1)^{n-1} - (B_{\diff}(\lambda) + 2 \epsilon) q^{n-1} > 1$.

Let $A$ be an $n$-dimensional abelian variety over $\F_q$.
Then $f_A(x) = x^n \, G(x+q/x)$ for some $G(x) = x^n + c_1 x^{n-1} + \cdots + c_n \in \Z[x]$  with all roots in $[-2q^{1/2},2q^{1/2}]$.
We have $c_i = O(q^{i/2})$ and $(a,b) \colonequals (q^{-1/2} c_1,q^{-1} c_2) \in \calS$.
Now 
\begin{equation}
\label{E:order of A}
   \#A(\F_q) = f_A(1) = G(q+1) = (q+1)^n + c_1 (q+1)^{n-1} + b q^{n-1} + O(q^{n-3/2}).
\end{equation}
Since $b = O(1)$, if $\#A(\F_q)$ lies in the interval \eqref{E:unrealizable for large q},
then $c_1 = r + O(1)$, so $a = q^{-1/2} c_1 = \lambda + O(q^{-1/2})$.
Then 
\[
      b \in [B_{\min}(a),B_{\max}(a)] \subset [B_{\min}(\lambda)-\epsilon/2,B_{\max}(\lambda)+\epsilon/2]
\]
by continuity, if $q$ is large enough.
If $c_1 \le r$, then the right side of \eqref{E:order of A} is too small to lie in \eqref{E:unrealizable for large q}.
If $c_1 \ge r+1$, then it is too large.
\end{proof}

\begin{lemma}
\label{L:large interval for large q}
Suppose that $n \ge 3$ and $\lambda \in \R$ satisfies $0 < \lambda < \lambda_1$.
Then for sufficiently large $q$, every integer in
\begin{equation}
\label{E:almost interval for large q}
   \bigl[\; q^n - \lambda q^{n-1/2} \;,\;  q^n + \lambda q^{n-1/2} \;\bigr]
\end{equation}
is $\#A(\F_q)$ for some $n$-dimensional abelian variety $A$ over $\F_q$.
\end{lemma}

\begin{proof}
By Lemma~\ref{L:B_diff}\eqref{I:B_diff}, $B_{\diff} >1$ on $[-\lambda,\lambda]$.
Choose $\epsilon>0$ so that $B_{\diff} > 1 + 2\epsilon$ on $[-\lambda,\lambda]$.
Let 
\[   
   \calS_\epsilon = \{\, (a,b) \in [-2n,2n] \times \R : B_{\min}(a) + \epsilon \le b \le B_{\max}(a) - \epsilon \, \}.
\]
Then $\calS_\epsilon$ is a compact subset of the interior of $\calS$.
Let $\calG_\epsilon$ be the inverse image of $\calS_\epsilon$ under $\calG_0 \surjects \calS$.
By Lemma~\ref{L:B_diff}\eqref{I:g_{a,b}}, 
$\calG_\epsilon$ is compact and consists of polynomials with distinct real roots in $(-2,2)$, 
so we can choose $\delta>0$ such that
any polynomial whose coefficients are within $\delta$ of some $g \in \calG_\epsilon$
again has distinct real roots in $(-2,2)$.

Suppose that $m$ is an integer in $[\; q^n - \lambda q^{n-1/2} \;,\;  q^n + \lambda q^{n-1/2} \;]$.
The rest of the proof relies on the following construction.

\begin{construction}\hfill
\begin{enumerate}[\upshape 1.]
\item \label{step:a} Let $a \in [-\lambda,\lambda]$ be such that $m = q^n + a q^{n-1/2}$.
\item \label{step:c_1 and b} Write $m = (q+1)^n + (c_1 + b) (q+1)^{n-1}$ with $c_1 \in \Z$ and $b \in [ B_{\min}(a) + \epsilon, B_{\max}(a) -\epsilon]$ (possible since $[ B_{\min}(a) + \epsilon, B_{\max}(a) -\epsilon]$ has length $>1$).   Then $(a,b) \in \calS_\epsilon$.
\item Choose $g \in \calG_\epsilon$ mapping to $(a,b)$.
\item Let $G(x) = q^{n/2} \, g(q^{-1/2} x) = x^n + q^{1/2} a x^{n-1} + q b x^{n-2} + \cdots \in \R[x]$.
\item Let $G_1$ be the same as $G$ except with the coefficient of $x^{n-1}$ changed to $c_1$.
\item For $i=2,\ldots,n$, let $G_i$ be the same as $G_{i-1}$ except with the coefficient of $x^{n-i}$
changed to the integer $c_i$ that makes $G_i(q+1)-m \in [0,(q+1)^{n-i})$.
\item \label{step:ordinary} Let $G_{\final} = G_n + s(x-(q+1))$, where $s \in \{0,1\}$ is chosen so that $p \nmid G_{\final}^{[0]}$.
\end{enumerate}
\end{construction}

\noindent\emph{Continuation of proof of Lemma~\ref{L:large interval for large q}.} 
We now bound the coefficients of $G_{\final} - G$ in order to prove that for $q$ large enough, the roots of $G_{\final}$ are still distinct and all in $[-2q^{1/2},2q^{1/2}]$.
Since $\calS_\epsilon$ is compact, $b$ is $O(1)$.
By Steps \ref{step:a} and~\ref{step:c_1 and b}, 
\begin{align}
\notag  q^n + a q^{n-1/2} &= m = (q+1)^n + (c_1+b)(q+1)^{n-1} = q^n + c_1 q^{n-1} + O(q^{n-1}),   \\
\label{E:c_1}   c_1 &= q^{1/2} a + O(1).
\end{align}
Now 
\begin{align*}
   G_1(q+1) &= (q+1)^n + c_1 (q+1)^{n-1} + qb (q+1)^{n-2} + O(q^{3/2}) (q+1)^{n-3} + \cdots + O(q^{n/2}) 1 \\
    &= (q+1)^n + (c_1 + b)(q+1)^{n-1} + O(q^{n-3/2}) \\
    &= m + O(q^{n-3/2}),
\end{align*}
so 
\begin{equation}
\label{E:c_2}
   c_2 - G^{[n-2]} = O(q^{n-3/2}) / (q+1)^{n-2} = O(q^{1/2}).
\end{equation}
Similarly, for $i=3,\ldots,n$, we have 
\begin{equation}
\label{E:c_3 etc}
   c_i - G^{[n-i]} = O((q+1)^{n-(i-1)})/(q+1)^{n-i} = O(q).
\end{equation}
Equations \eqref{E:c_1}, \eqref{E:c_2}, and \eqref{E:c_3 etc} imply that
\[
    G_n^{[n-i]} - G^{[n-i]} = O(q^{(i-1)/2})
\]
for all $i \ge 1$.
Since $n \ge 3$, the same holds with $G_n$ replaced by $G_{\final}$.
Thus the coefficients of $g_{\final}(x) = q^{-n/2} \, G_{\final}(q^{1/2} x)$
are within $O(q^{-1/2}) < \delta$ of the corresponding coefficients of $g$ 
if $q$ is sufficiently large, so $g_{\final}$ has all its roots in $[-2,2]$.
Thus $G_{\final}$ has all its roots in $[-2q^{1/2},2q^{1/2}]$.
By construction, $G_{\final} \in \Z[x]$.
Also $G_{\final}(q+1) - m = G_n(q+1) - m \in [0,1)$, so $G_{\final}(q+1)=m$.

Let $f(x)=x^n \, G_{\final}(x+q/x) \in \Z[x]$.
We have $f^{[n]} \equiv G_{\final}^{[0]} \not\equiv 0 \pmod{p}$.
By Theorem~\ref{T:Honda-Tate}, $f=f_A$ for some $n$-dimensional ordinary abelian variety over $\F_q$.
Finally, $\#A(\F_q)=f(1)=G_{\final}(q+1)=m$.
\end{proof}

\begin{proof}[Proof of Theorem~\ref{T:large q}]
Lemma~\ref{L:large interval for large q} shows that all integers in
$[ q^n - \lambda q^{n-1/2} , q^n + \lambda q^{n-1/2} ]$
are realizable for $\lambda$ that can approach $\lambda_1$ from below as $q \to \infty$.
Lemma~\ref{L:unrealizable} shows, on the other hand, 
that for any $\mu$ with $|\mu| > \lambda_1$, 
there are unrealizable integers within $O(q^{n-1})$ of $(q+1)^n + \mu q^{n-1/2}$ 
if $q$ is sufficiently large.
These imply Theorem~\ref{T:large q}.
\end{proof}

\begin{remark} \label{R:n=2}
Suppose $n=2$.
Theorem~\ref{T:large q} holds without change if $q$ tends to $\infty$ through primes only: 
the proof of Lemma~\ref{L:large interval for large q} works 
if we omit Step~\ref{step:ordinary}, because of the last sentence of Remark~\ref{R:Honda-Tate}.

On the other hand, if $q$ tends to $\infty$ through non-prime prime powers,
then Theorem~\ref{T:large q} holds with $\lambda_1$ replaced by the smaller value 
$\lambda_2 = 4-2\sqrt{2}$, as we now explain.
In Lemma~\ref{L:large interval for large q}, if $0 < \lambda < \lambda_2$,
then $B_{\diff} > 2$ on $[-\lambda, \lambda]$,
so there are at least \emph{two} consecutive integer possibilities for $c_1$,
and at least one of them will lead to a polynomial $f$ 
for which \eqref{I:middle coefficient} in Theorem~\ref{T:Honda-Tate} holds.
Meanwhile, in Lemma~\ref{L:unrealizable}, if $\lambda_2 < |\mu| < 2n$,
so that $B_{\diff}(\mu) < 2$, then there exists $\epsilon>0$
such that if $q$ is sufficiently large, 
and $r$ is the multiple of $p$ nearest $\mu q^{1/2}$,
then any integer of the form $m=(q+1)^2 + r (q+1) +c_2$ in
\[
\bigl[ \; (q+1)^2 + (r-1) (q+1) + (B_{\max}(\mu) + \epsilon) q \; , \; (q+1)^2 + (r+1) (q+1) + (B_{\min}(\mu) - \epsilon) q \; \bigr]
\]
with $p \mid c_2$ and $p^2 \nmid c_2$
is not $\#A(\F_q)$ for any abelian surface $A$ over $\F_q$,
because the only monic quadratic polynomial $G(x) \in \Z[x]$ with roots in $[-2q^{1/2},2q^{1/2}]$
such that $G(q+1)=m$ is $x^2+rx+c_2$, which is Eisenstein at $p$,
which implies that the polynomial $f(x) \colonequals x^2 \, G(x+q/x)$ fails condition~(d$'$) in Remark~\ref{R:Honda-Tate}.
\end{remark}

%****************************************************************************
\section{Effective bounds} 
\label{S:effective}

Given $q$ and $n$, we have given three ways to construct polynomials
that realize a large interval of integers as $\#A(\F_q)$ for $A$ of dimension $n$ over $\F_q$:
\begin{itemize}
\item Section~\ref{S:simple proof} gave a quick construction 
that realized intervals wide enough to cover all sufficiently large integers
as $n$ varies, but not wide enough to be asymptotically close to optimal.
\item Section~\ref{sec:range} gave a more subtle construction that
gave a much wider interval, but it is too complicated to analyze
explicitly to make all the big-$O$ constants explicit.
\item Section~\ref{S:large q} gave a method that again is asymptotically good,
but only when $q$ is large compared to $n$.
\end{itemize}
In this section, we present a \emph{fourth} construction that,
while not asymptotically as good as the construction of Section~\ref{sec:range},
realizes a wide interval for any $q$ and sufficiently large $n$,
and is still simple enough to analyze fully.

Given $q$, $n \ge 2$, and an integer $m$ in $[q^{n-1/2},q^{n+1/2})$,
the plan is to find a power series $j(z) \in z\R[[z]]$ such that
$j(1/q) = \log(m/q^n)$ and $\exp(j(z)) \in \Z[[z]]$;
then we truncate $\exp(j(z))$ to a degree~$n$ polynomial 
and adjust the coefficients of $x^{n-1}$ and $x^n$
to produce a polynomial $h(z)$ such that $\widehat{h}(1)=m$ and $p \nmid \widehat{h}^{[n]}$.
This should work well, since $\exp(j(z))$ is automatically nonvanishing on $D$, 
and if its coefficients are not too large, then the nonvanishing should persist after truncating and adjusting.

\begin{construction}
\label{C:exponential}
\hfill
\begin{enumerate}[\upshape 1.]
\item \label{step:integer coefficient} 
For $i=1,2,\ldots,n-1$, let $c_i$ be the real number such that 
\[
   \log(m/q^n) - c_1 q^{-1} - \cdots - c_i q^{-i} \in [-q^{-i}/2,q^{-i}/2)
\]
and such that the coefficient of $z^i$ in $\exp(c_1 z + \cdots + c_i z^i)$ is an integer;
for the existence and uniqueness of $c_i$, see the proof of Lemma~\ref{L:c_1 is at most s}.
\item \label{step:c_n}
Let $c_n \in \R$ be such that $\log(m/q^n) - c_1 q^{-1} - \cdots - c_n q^{-n} = 0$.
\item Let $h_0(z) \in \R[z]$ be the degree~$n$ Taylor polynomial of $\exp(c_1 z + \cdots + c_n z^n)$.
\item Let $h_1(z) = h_0(z) + k z^n/2$, where $k \in \R$ is chosen to make $\widehat{h}_1(1)=m$.
\item \label{step:ordinary adjustment}
Let $h$ be $h_1$ or $h_1 + z^{n-1} - ((q+1)/2) z^n$, whichever makes $p \nmid \widehat{h}^{[n]}$.
\item Let $A$ be an abelian variety with $f_A=\widehat{h}$, if one exists.  (If $h$ is nonvanishing on $D$, then such an $A$ is guaranteed to exist and $\widehat{h}$ is squarefree by Remark~\ref{R:squarefree}.)
% \item If $h$ is nonvanishing on $D$, let $A$ be such that $f_A = \widehat{h}$. 
\end{enumerate}
\end{construction}

Let $s = \lfloor \frac 12 q \log q + \frac12 \rfloor$.

\begin{lemma}
\label{L:c_1 is at most s}
We have $|c_1| \le s$ and $|c_i| \le (q+1)/2$ for $i=2,\ldots,n$.
\end{lemma}

\begin{proof}
Since $m \in [q^{n-1/2},q^{n+1/2})$, we have $\log(m/q^n) \in [-\frac12 \log q, \frac12 \log q)$,
and Step~\ref{step:integer coefficient} says that $c_1$ is the integer in the interval $q \log(m/q^n) + [-\frac12,\frac12)$,
so $|c_1| \le s$.

For $i \le n-1$, let $\epsilon_i =\log(m/q^n) - c_1 q^{-1} - \cdots - c_i q^{-i}$,
so $\epsilon_i = \epsilon_{i-1} - c_i q^{-i}$; then
$\epsilon_{i-1} \in [-q^{-(i-1)}/2,q^{-(i-1)}/2)$,
so the condition $\epsilon_i \in [-q^{-i}/2,q^{-i}/2)$  in Step~\ref{step:integer coefficient}
constrains $c_i$ to a half-open interval of length $1$ contained in $[-(q+1)/2,(q+1)/2]$,
while the integer coefficient condition in Step~\ref{step:integer coefficient} constrains $c_i$ to a coset of $\Z$ in $\R$; thus a unique $c_i$ exists, and $|c_i| \le (q+1)/2$.
Finally, $c_n = q^n \epsilon_{n-1} \in q^n[-q^{-(n-1)}/2,q^{-(n-1)}/2) = [-q/2,q/2)$.
\end{proof}

Let $j(z) = c_1 z + \cdots + c_n z^n$. 
To bound the difference between $\exp(j(z))$ and its degree~$n$ Taylor polynomial, 
we consider the worst case:  let 
\[
   J(z) \colonequals \exp\left( sz + \frac{q+1}{2} \frac{z^2}{1-z} \right) = J_{\le n}(z) + J_{>n}(z) \quad \in \R_{\ge 0}[[z]],
\]
where $J_{\le n}$ is the degree~$n$ Taylor polynomial, 
and $J_{>n}$ is the remainder power series consisting of terms of degree~$>n$.
By Lemma~\ref{L:c_1 is at most s}, $|(\exp j(z))^{[i]}| \le J^{[i]}$.

\begin{proposition}
\label{P:J inequality}
Let $q$ be a prime power.  For $n \ge 2$ and $m \in [q^{n-1/2},q^{n+1/2})$, if 
\begin{equation}
\label{E:J inequality}
    J_{>n}(q^{-1/2}) + \frac{q^{n/2}}{2} J_{>n}(q^{-1}) + \frac{q^{-n/2}}{2} J_{\le n}(1) +  \frac{(q^{1/2}+1)^2}{2} q^{-n/2} < \frac{1}{J(q^{-1/2})},
\end{equation}
then Construction~\ref{C:exponential} produces an ordinary $n$-dimensional $A$ over $\F_q$ with $\#A(\F_q)=m$.
\end{proposition}

\begin{proof}
By Step~\ref{step:c_n}, $\exp j(q^{-1}) = m/q^n$, so
\begin{gather*}
   \lvert m - q^n h_0(q^{-1}) \rvert \;=\;  q^n \lvert \exp j(q^{-1}) - h_0(q^{-1}) \rvert \;\le\; q^n \, J_{>n}(q^{-1}) \\
   \lvert h_0(1) \rvert \;\le\; J_{\le n}(1) \\
   |k| \;=\; \lvert \widehat{h}(1)- \widehat{h}_0(1) \rvert \;\le\; \lvert m - q^n h_0(q^{-1}) - h_0(1) \rvert \;\le\; q^n \, J_{>n}(q^{-1}) + J_{\le n}(1).
\end{gather*}
On $D$, 
\begin{align*}
   \lvert \exp j(z) \rvert &= \exp(\re j(z)) \ge \exp\left(-s q^{-1/2} - \frac{q+1}{2} q^{-2/2} - \cdots - \frac{q+1}{2} q^{-n/2} \right) \ge \frac{1}{J(q^{-1/2})} \\
   \lvert h_0(z) \rvert &\ge \lvert \exp j(z) \rvert - J_{>n}(q^{-1/2}) \\
      \lvert h_1(z) \rvert &\ge \lvert h_0(z) \rvert - \frac{k}{2} q^{-n/2} \\
      \lvert h(z) \rvert &\ge \lvert h_1(z) \rvert - q^{-(n-1)/2} - \frac{q+1}{2} q^{-n/2}
           = \lvert h_1(z) \rvert - \frac{(q^{1/2}+1)^2}{2} q^{-n/2}.
\end{align*}
Combining the previous five inequalities yields
\begin{align*}
   |h(z)| &\ge \frac{1}{J(q^{-1/2})} - J_{>n}(q^{-1/2}) - \frac{q^{n/2}}{2} J_{>n}(q^{-1})  - \frac{q^{-n/2}}{2} J_{\le n}(1) - q^{-(n-1)/2} - \frac{q+1}{2} q^{-n/2},
\end{align*}
so \eqref{E:J inequality} implies that $h$ is nonvanishing on $D$.
Theorem~\ref{T:Honda-Tate} produces $A$.
\end{proof}
The following weaker statement has the advantage that if any hypothesis holds for one $n$,
it clearly holds for all larger $n$:
\begin{corollary}
\label{C:weakenings}
Let $q$ be a prime power.  For $n \ge 2$ and $m \in [q^{n-1/2},q^{n+1/2})$, if any of 
\begin{gather}
\label{E:simpler J inequality}
    (1+q^{-1/2}/2) J_{>n}(q^{-1/2}) + \tfrac{1}{2} \left( \tfrac{4}{3} q^{-1/2} \right)^n J(\tfrac{3}{4}) +  \frac{(q^{1/2}+1)^2}{2} q^{-n/2} < \frac{1}{J(q^{-1/2})}, \\
\label{E:q at least 7}
  q \ge 7 \quad \textup{and} \quad  2^{n-1} >  J(q^{-1/2}) J(2q^{-1/2}), \quad \textup{or}\\
\label{E:q at least 16}
  q \ge 16 \quad \textup{and} \quad n > 3 q^{1/2} \log q - 1/2,
\end{gather}
hold, then Construction~\ref{C:exponential} produces an ordinary $n$-dimensional $A$ over $\F_q$ with $\#A(\F_q)=m$.
\end{corollary}

\begin{proof}
In \eqref{E:J inequality}, $J_{>n}(q^{-1})  \le q^{-(n+1)/2} J_{>n}(q^{-1/2})$ 
(this holds termwise for any power series $J$ with nonnegative coefficients).
Similarly $J_{\le n}(1) \le (\frac{4}{3})^n J_{\le n}(3/4) \le (\frac{4}{3})^n J(3/4)$.
Hence the left side of \eqref{E:J inequality} is at most the left side of \eqref{E:simpler J inequality}.
Thus, if \eqref{E:simpler J inequality} holds, Proposition~\ref{P:J inequality} applies.

Now suppose that $q \ge 7$ and $2^{n-1} > J(q^{-1/2} )J(2 q^{-1/2})$.
First,
\[
   2^{n-1} > J(q^{-1/2}) J(2q^{-1/2}) \ge \exp(3sq^{-1/2}) \ge \exp(3 q^{-1/2} (q \log q -1)/2 ) \ge 2^9,
\]
so $n \ge 10$.
We use 
\begin{equation}
\label{E:four inequalities}
\begin{aligned}
   J_{>n}(q^{-1/2}) &\le 2^{-(n+1)} J(2q^{-1/2}) \\
   J_{>n}(q^{-1}) &\le (2 q^{1/2})^{-(n+1)} J(2q^{-1/2}) \\
   J_{\le n}(1) &\le 1 + (q^{1/2}/2)^n ( J(2q^{-1/2}) -1) \\
   (q^{1/2}+1)^2 &\le (q^{1/2}/2)^n -1;
\end{aligned}
\end{equation}
the first three are proved termwise, and the last follows from 
the inequality $(2u+1)^2 \le u^{10}-1$ for $u \colonequals q^{1/2}/2 \ge 7^{1/2}/2$. 
By \eqref{E:four inequalities}, the left side of \eqref{E:J inequality} is at most
\begin{align*}
   \lefteqn{2^{-(n+1)} J(2q^{-1/2}) + \frac{q^{n/2}}{2} (2q^{1/2})^{-(n+1)} J(2q^{-1/2})} & \\
   & \qquad {} + \frac{q^{-n/2}}{2} \Bigl( (q^{1/2}/2)^n J(2q^{-1/2}) + 1 - (q^{1/2}/2)^n \Bigr) + \frac{q^{-n/2}}{2} \Bigl( (q^{1/2}/2)^n - 1 \Bigr) \\
   &= 2^{-(n+1)}(2 + q^{-1/2}/2) \, J(2q^{-1/2}) \\
   &\le 2^{1-n} J(2 q^{-1/2}) \\
   &< \frac{1}{J(q^{-1/2})},
\end{align*}
by hypothesis, so Proposition~\ref{P:J inequality} applies.

Finally, suppose that $q \ge 16$ and $n > 3 q^{1/2} \log q - 1/2$.
Then
\begin{align}
\notag     s &\le (q \log q + 1)/2 \\
\notag     \log\bigl(J(q^{-1/2}) J(2q^{-1/2}) \bigr) &\le 3 \Bigl( \frac{q \log q +1}{2} \Bigr) q^{-1/2} + \frac{q+1}{2} \Bigl(\frac{q^{-1}}{1-q^{-1/2}} + \frac{4 q^{-1}}{1-2q^{-1/2}} \Bigr) \\
\label{E:log inequality}    &\le (3 q^{1/2} \log q - 3/2) \log 2 \\
\notag     &< (n-1) \log 2,
\end{align}
so \eqref{E:q at least 7} holds;
to prove \eqref{E:log inequality}, check numerically for $16 \le q \le 100$
and for $q>100$ use
\begin{align*}
   \lefteqn{\frac{3}{2} q^{-1/2} + \frac{q+1}{2} \Bigl(\frac{q^{-1}}{1-q^{-1/2}} + \frac{4 q^{-1}}{1-2q^{-1/2}} \Bigr) + \frac{3}{2} \log 2} & \\
   &\le \frac{3}{2}(0.1) + q \Bigl(\frac{q^{-1}}{0.9} + \frac{4 q^{-1}}{0.8} \Bigr) + \frac{3}{2} \log 2 < 8 < (3 \log 2 - 3/2) q^{1/2} \log q.\qedhere
\end{align*}
\end{proof}

Corollary~\ref{C:weakenings} proves Theorem~\ref{T:effective}\eqref{I:constant 3} for $q \ge 16$.
Also, for each $q<16$ it provides an $n$ such that all integers $\ge q^{n-1/2}$ are realizable,
but too many integers remain to be checked one at a time.
Therefore we describe a construction allowing us to realize larger intervals of integers all at once.
The plan is to start with $h$ such that $\widehat{h} = f_A$ for some $A$ with $\#A(\F_q)=m$,
and then to replace $h$ by $h + \sum_{i=r}^n c_i z^i$ for some $r$ and small $c_i$
(and then adjust to make $p \nmid \widehat{h}^{[n]}$ again);
as the $c_i$ vary, we realize all integers in an interval.

\begin{construction}
\label{C:interval}
Suppose that we are given $q$, $n$, $m$, and a polynomial $h\in 1 + x \Z[x]$ of degree $<2n$ with $\widehat{h}(1)=m$
(given by Construction~\ref{C:exponential} or otherwise).
\begin{enumerate}[\upshape 1.]
\item \label{step:zeros of h}
Compute the complex zeros of $h$ and check that none of them are in $D$.
(More precisely: Compute small balls containing the zeros, and check that none of them intersect $D$.)
\item
Compute the complex zeros $\alpha$ of the derivative of $h(z) h(1/(qz))$,
evaluate $|h|$ at each $\alpha$ on the boundary $\del D$, and let $\mu$ be the minimum of these values; see the proof of Lemma~\ref{L:interval construction}.
(More precisely: Compute small balls around these zeros, and let $\mu$ be a lower bound for $|h|$ on all these balls that intersect $\del D$; if $h=1$, then let $\mu=1$.)
\item \label{step:mu_ord}
Let $\mu_{\ord} = \mu - q^{-(n-1)/2} - ((q+1)/2) q^{-n/2}$; abort if $\mu_{\ord} \le 0$.
\item Let $r$ be the smallest positive integer $\le n+1$ such that $\sum_{i=r}^n \lfloor q/2 \rfloor q^{-i/2} < \mu_{\ord}$. 
\item
Let $N = \lfloor q/2 \rfloor \sum_{j=r}^n (q^{n-j}+1) = \lfloor q/2 \rfloor \left( \frac{q^{n-r+1}-1}{q-1} + (n-r+1) \right)$.
\item \label{step:output interval}
Return the interval $\bigl[\,\widehat{h}(1) - N ,\widehat{h}(1) + N\bigr]$.
\end{enumerate}
\end{construction}

\begin{lemma}
\label{L:interval construction}
In Construction~\ref{C:interval}, if Steps \ref{step:zeros of h} and~\ref{step:mu_ord} succeed,
then every integer in the interval of Step~\ref{step:output interval} 
is $\#A(\F_q)$ for some ordinary abelian variety of dimension $n$ over $\F_q$.
\end{lemma}

\begin{proof}
Suppose that $h$ has no zeros in $D$.
Then the minimum of $|h|$ occurs on $\del D$,
where $|h|^2 = h(z) \, h(1/(qz))$,
and this minimum occurs at a point where the derivative of $h(z) \, h(1/(qz))$ is $0$.
Thus $|h| \ge \mu$ on $D$.

Suppose that $H = h + \sum_{i=r}^n c_i z^i$ where $|c_i| \le q/2$ for all $i$, and $c_i \in \Z$ for all $i$ except $n$, and $c_n \in \frac12\Z$.
The choice of $r$ guarantees that $|H| > 0$ on $D$, even if we add $z^{n-1} - ((q+1)/2) z^n$ to $H$
if necessary to make $p \nmid \widehat{H}^{[n]}$, 
so $\widehat{H}(1)$ is realizable.
To realize an integer $\widehat{h}(1) + M$ with $|M| \leq N$, write $M$ as $\sum_{j=r}^n c_j(q^{n-j} + 1)$ with $|c_j| \leq  \lfloor q/2 \rfloor$, $c_j \in \Z$ for all $j \neq n$, and $c_n \in \tfrac12 \Z$, by greedily choosing $c_r$, $c_{r+1}$, \dots, one at a time.
\end{proof}

\begin{proof}[Proof of Theorem~\ref{T:effective} and Remarks \ref{R:realizing ordinary}, \ref{R:7 and 8}, and~\ref{R:realizing squarefree}]
In this proof, given $q$, a positive integer is called realizable
if it equals $\#A(\F_q)$ for some ordinary abelian variety $A$ over $\F_q$ with $f_A$ squarefree.
The case $q=2$ is done by \cite{Howe-Kedlaya-preprint}.
Criterion~\eqref{E:q at least 16} of Corollary~\ref{C:weakenings} proves Theorem~\ref{T:effective}\eqref{I:constant 3} for $q \ge 16$.
For each $q<16$, we numerically find $n \ge 2$ such that \eqref{E:simpler J inequality} holds;
then we check smaller values of $n$ to find the smallest $n_0$ such that \eqref{E:J inequality} holds for all $n \ge n_0$.
(It turns out that $n_0 \le 25$ for each $q<16$.)
For $q \in \{11,13\}$, it turns out that $q^{3 \sqrt{q} \log q} > q^{n_0-1/2}$, which proves Theorem~\ref{T:effective}\eqref{I:constant 3} for these $q$.

For $3 \le q \le 9$, we use variants of Construction~\ref{C:exponential} and \ref{C:interval} 
to realize all integers in an interval $[M_q,q^{n_0-1/2}]$.
For $q \in \{8,9\}$, we have $M_q \le q^{3 \sqrt{q} \log q}$, 
which proves Theorem~\ref{T:effective}\eqref{I:constant 3} for these $q$.
For $q \in \{3,4,5,7\}$, we 
use the algorithm of~\cite{Kedlaya-2008} (implemented at \url{https://github.com/kedlaya/root-unitary})
to exhaust over the polynomials $f_A$ for abelian varieties $A$ of dimension $\le 4$
to realize all integers $<M_q$ 
with the exception of those listed in Remarks \ref{R:realizing ordinary}, \ref{R:7 and 8}, and~\ref{R:realizing squarefree}.
Neither are these exceptions realized by abelian varieties of dimension $\ge 5$,
because they are all less than $(\sqrt{q}-1)^{10}$.
The calculations in this paragraph took 7.19 CPU hours on an Intel Core i7-9750H CPU @ 2.60GHz.
See \url{https://github.com/edgarcosta/abvar-fq-orders} for the code and further details.
\end{proof}

Some calculations were checked against the database of isogeny classes of abelian varieties over finite fields in the L-functions and Modular Forms Database \cite{abvardb, LMFDB}. 

%***************************************************************************t

\appendix   
\section{Optimality of a potential function} 
\label{A:potential theory}

\subsection{Polynomials}
\label{S:polynomials}

The goal of this appendix is to prove the following.

\begin{proposition}
\label{P:M(r)}
Choose $c$ in the interval $(0, 1)$. 
For $d \ge 1$, 
let $\FF(d, c)$ be the set of complex polynomials $f$ of degree $d$ satisfying
$f(0)=1$ and $|f(w)|^{1/d} \ge c$ for all $w \in \C_{\le 1}$.
On $(-\infty,1]$ define the decreasing continuous function
\[
    M(r) \colonequals \frac{1 - r + \sqrt{(1 - r)^2 + 4rc^2}}{2}.
\]
\begin{enumerate}[\upshape (a)]
\item \label{I:F(n,c) bounds}
For any $f \in \FF(d, c)$, we have  
\begin{align}
\label{E:polynomial lower bound}
    |f(r)|^{1/d} &\ge M(r) \quad\;\;\:\textup{ for all $r \in [0,1]$,} \\
\label{E:polynomial upper bound}
    |f(r)|^{1/d} &\le M(-r) \quad\textup{ for all $r \in [0,\infty)$.} 
\end{align}
\item \label{I:F(n,c) limits}
There exist polynomials $f_1,f_2,\ldots$ with $f_d \in \FF(d,c)$ such that for every $r \in (-\infty,1]$,
\begin{equation}
\label{E:inf sup}
   \lim_{d \to \infty} |f_d(r)|^{1/d} = M(r).
\end{equation}
$($Thus \eqref{E:polynomial lower bound} is asymptotically sharp, and \eqref{E:polynomial upper bound} is too
since $f_d(-z) \in \FF(d,c)$.$)$
\end{enumerate}
\end{proposition}

\begin{remark}
For $r>1$, the lower bound in \eqref{E:polynomial lower bound} is simply $0$ since the function $f(z) \colonequals 1-(z/r)^d$ is in $\FF(d,c)$.
\end{remark}

\begin{remark}
If $f$ is in $\FF(d,c)$, then so is $f(uz)$ for any $u \in \C$ with $|u|=1$.
Thus Proposition~\ref{P:M(r)} implies the same results for $f(w)$
for any complex number $w$ satisfying $|w|=r$.
\end{remark}

Outside the trivial case $r = 0$ and the case $r = 1$, which was handled in detail in \cite{RuVa86}, 
Proposition~\ref{P:M(r)} appears to be new.

\begin{remark}
Choose a prime power $q$, and take $c = q^{-1/4}$ and $r = q^{-1/2}$. 
Let $I_1, I_2, \dots$ be an increasing sequence of closed intervals with union $\Ia$. For each positive integer $k$, let $P_k$ be a polynomial constructed as in Proposition~\ref{P:construction of P} from the interval $I_k$. Then the polynomials $P_k(q^{1/2}z) \in \FF(\text{deg} P_k, c)$ have limits as in \eqref{E:inf sup}.
In the other direction, Proposition~\ref{P:M(r)}\eqref{I:F(n,c) bounds} shows that we could not hope to construct polynomials satisfying the conditions of Proposition~\ref{P:construction of P} for intervals larger than $\Ia$.
\end{remark}

\subsection{Potential functions}
\label{S:potential functions}

Given a nonconstant polynomial $f$, let $\mu$ be the uniform probability measure on the set of zeros of $f$, counted with multiplicity.
Then $\log |f(z)|^{1/d}$ equals $\int \log|w-z|\,d\mu(w)$ minus a constant,
so Proposition~\ref{P:M(r)} can be reformulated in terms of $\mu$.
This suggests a generalization in which $\mu$ is allowed to be any compactly supported probability measure
on $\C_{\ge 1}$.
In fact, this generalization, formalized as Proposition~\ref{P:potential version} below, 
is \emph{equivalent} to Proposition~\ref{P:M(r)}.

\begin{definition}[{\cite[I.1]{Saff-Totik-1997}}]
Let $\Sigma$ be a compact subset of $\C$.
Let $\mathcal{M}(\Sigma)$ be the set of (Borel) probability measures on $\C$ with support contained in $\Sigma$. 
For $\mu \in \mathcal{M}(\Sigma)$, define the \defi{potential function} $U^{\mu} \colon \C \to \R \union \{\infty\}$ by
\[
   U^{\mu}(z) \colonequals \int_{\C} -\log|w - z|\, d\mu(w).
\]
\end{definition}

For a polynomial $F$ with nonzero constant term, the polynomial $f(z) \colonequals F(z)/F(0)$ satisfies $f(0)=1$, as required in the definition of $\FF(d,c)$.  Analogously, we will consider $U^{\mu}(z)-U^\mu(0)$.

\begin{proposition}
\label{P:potential version}
Choose $c$ in $(0, 1)$, and let $\MM(c)$ be the set of probability measures $\mu$ with compact support contained in $\C_{\ge 1}$ such that 
\begin{equation}
    \label{E:c_assumption}
    U^{\mu}(z) - U^\mu(0) \;\le\; - \log c \qquad\text{for all $z \in \C_{\le 1}$.}
\end{equation}
\begin{enumerate}[\upshape (a)]
\item \label{I:bounds on U^mu}
For any $\mu \in \MM(c)$,  
\begin{align}
\label{E:potential upper bound}
    U^{\mu}(r) - U^{\mu}(0) &\le -\log M(r) \quad\;\;\text{ for all $r \in [0,1]$,} \\
\label{E:potential lower bound}
    U^{\mu}(r) - U^{\mu}(0) &\ge -\log M(-r) \quad\text{ for all $r \in [0,\infty)$.} 
\end{align}
\item \label{I:extreme measures}
Let $\Sigma_c$ be the arc $\{z \in \C : |z|=1 \textup{ and } |z-1| \le 2c\}$.
There exists a measure $\mu_c \in \MM(c)$ supported on $\Sigma_c$ such that for every $r \in (-\infty,1]$,
\[
   U^{\mu_c}(r) - U^{\mu_c}(0) = -\log M(r). 
\]
\end{enumerate}
\end{proposition}

\begin{proof}[Proof that Proposition~\ref{P:potential version} implies Proposition~\ref{P:M(r)}]
If Proposition~\ref{P:potential version}\eqref{I:bounds on U^mu} holds, 
apply it to the uniform probability measure $\mu$
on the zeros of $f \in \FF(d,c)$ counted with multiplicity, 
and apply $x \mapsto e^{-x}$ to \eqref{E:potential upper bound} and~\eqref{E:potential lower bound} to get Proposition~\ref{P:M(r)}\eqref{I:F(n,c) bounds}.

Now suppose in addition that Proposition~\ref{P:potential version}\eqref{I:extreme measures} holds.
Fix $r \in [0,1]$.
For $\lambda \in (0,1)$, Proposition~\ref{P:potential version}\eqref{I:bounds on U^mu} (and rotational symmetry) shows that
for $z \in \C_{\le \lambda}$, 
\begin{equation}
\label{E:potential lower bound on disk}
    U^{\mu_c}(z)-U^{\mu_c}(0) \le -\log M(\lambda),
\end{equation}
which is \emph{strictly} less than $-\log c$.
By approximating $\mu_c$ by uniform probability measures supported on finite subsets of $\C_{\ge 1}$,
we find a sequence of polynomials $p_1, p_2, \dots$ such that 
\begin{itemize}
\item $p_d$ has degree $d$, has all roots in $\C_{\ge 1}$, and satisfies $p_d(0)=1$; and
\item on each compact subset of $\C \setminus \Sigma_c$, the sequence $-\log |p_d(z)|^{1/d}$ converges uniformly
to $U^{\mu_c}(z)-U^{\mu_c}(0)$.
\end{itemize}
By \eqref{E:potential lower bound on disk} and uniform convergence, 
for any $\lambda<1$, if $d$ is sufficiently large, 
then $-\log |p_d(z)|^{1/d} \le -\log c$ on $\C_{\le \lambda}$, 
so the polynomial $f_d(z) \colonequals p_d(\lambda z)$ lies in $\FF(d,c)$.
Then $|f_d(r)|^{1/d} \to \exp(-(U^{\mu_c}(\lambda r)-U^{\mu_c}(0))) = M(\lambda r)$ 
uniformly on each compact subset of $(-\infty,1]$.
By repeating the argument for each $\lambda \in (0,1)$ to obtain $f_{d,\lambda}$,
and then letting $\lambda$ tend to $1$ sufficiently slowly with $d$,
we obtain polynomials satisfying Proposition~\ref{P:M(r)}\eqref{I:F(n,c) limits}.
\end{proof}

If $\mu$ is supported on the unit circle, then $U^\mu(0)=0$ by definition.
In proving Proposition~\ref{P:potential version}\eqref{I:bounds on U^mu}, 
the following lets us assume that $\mu$ is supported on the unit circle.

\begin{lemma}
\label{L:Balay}
Given a compactly supported probability measure $\mu$ on $\C_{\ge 1}$, there is a probability measure $\widehat{\mu}$ supported on the unit circle $\Sigma$ such that
\[U^{\mu}(z) - U^{\mu}(0) = U^{\widehat{\mu}}(z)\quad\text{whenever}\quad |z| < 1\] 
and
\[U^{\mu}(z) - U^{\mu}(0) \ge U^{\widehat{\mu}}(z)\quad\text{whenever}\quad |z| \ge 1.\]
\end{lemma}

\begin{proof}
Write $\mu$ as a sum of nonnegative measures $\mu_\Sigma + \mu'$, where $\mu_\Sigma$ is supported on the circle and $\mu'(\Sigma)=0$.
Apply ``balayage'' (\cite[Theorem~II.4.7]{Saff-Totik-1997}) to $\mu'$ to produce $\widehat{\mu'}$ supported on the circle,
and let $\widehat{\mu} = \mu_\Sigma + \widehat{\mu'}$.
\end{proof}

\subsection{Equilibrium measures}
\label{S:equilibrium measures}

\begin{definition}
\label{D:equilibrium}
Suppose that $\Sigma$ is of positive capacity, as defined in \cite[Definition 1.5]{Saff-Totik-1997}; 
this holds if $\Sigma$ contains a line segment or circular arc of positive length, for example. 
The \defi{energy} of $\mu \in \mathcal{M}(\Sigma)$ is $\int_{\Sigma} U^{\mu}(z) \, d\mu(z)$.
There is a unique energy-minimizing measure $\mu \in \mathcal{M}(\Sigma)$, 
called the \defi{equilibrium measure} on $\Sigma$.
More generally, for any continuous function $Q \colon \Sigma \to \R$,
there is a unique measure $\mu$ in $\mathcal{M}(\Sigma)$ minimizing the \defi{weighted energy}
\[
     E_Q(\mu) \colonequals \int_{\Sigma} \big(U^{\mu}(z) + 2 \, Q(z)\big)\, d\mu(z), 
\]
and this $\mu$ is called the \defi{weighted equilibrium measure} for $Q$ on $\Sigma$.
\end{definition}

From now on, $\Sigma$ denotes the unit circle, and $\kappa \colon \C^\times \to \C$ denotes the rational function $\kappa(z) \colonequals z+z^{-1}$, which maps $\Sigma$ onto the interval $[-2,2]$.

\begin{lemma}
\label{L:pushforward}
\hfill
\begin{enumerate}[\upshape (a)]
\item \label{I:bijection}
The map $\mu \mapsto \kap_* \mu$ sending measures to their pushforwards under $\kap$ is a bijection from 
the set of complex-conjugation-invariant probability measures on the unit circle $\Sigma$
to $\calM([-2,2])$.
\item \label{I:compare potentials}
For $\mu$ as in \eqref{I:bijection}, we have $U^{\kap_*\mu}(\kap(z)) = 2 \, U^\mu(z) + \log |z|$ for all $z \in \C^\times$.
\item \label{I:compare equilibrium measures}
Let $\Sigma'$ be a positive-capacity complex-conjugation-invariant compact subset of $\Sigma$.
Let $\alpha \in \R$ and $r \in \R^\times \setminus \Sigma'$.
Under $\kap_*$, the equilibrium measure on $\Sigma'$ for weight $Q(z) \colonequals \alpha \log|z - r|$
corresponds to the equilibrium measure on $\kap(\Sigma')$ for weight $R(z) \colonequals \alpha \log|z - \kap(r)|$.
\end{enumerate}
\end{lemma}

\begin{proof}\hfill
\begin{enumerate}[\upshape (a)]
\item The map $\kap$ induces an isomorphism from the $\sigma$-algebra of complex-conjugation-invariant Borel subsets of $\Sigma$ to the $\sigma$-algebra of Borel subsets of $[-2,2]$.
\item For $w,z \in \C^\times$ we have
\begin{equation}
    \label{E:smushing_circles}
    \kap(w)-\kap(z) = -\frac{1}{z} (w-z)(w^{-1}-z).
\end{equation}
The claim follows by applying $-\log |\;|$ and integrating against $d\mu(w)$.
\item 
Renaming variables in \eqref{E:smushing_circles} yields $R(\kap(z)) = 2 \, Q(z) - \alpha \log |r|$.
By symmetry, the only measures on $\Sigma'$ we need to consider are those that are complex-conjugation-invariant.
For such $\mu$,
\begin{align*}
   E_R(\kap_* \mu) &= \int_{\Sigma'} \big( U^{\kap_* \mu}(\kap(z)) + 2 \, R(\kap(z)) \big) \, d\mu(z) \\
   &= \int_{\Sigma'} \big( 2 \, U^{\mu}(z) + 4 \, Q(z) - 2 \alpha \log |r| \big) \, d\mu(z) \qquad\textup{(by \eqref{I:compare potentials})}\\
   &= 2 \, E_Q(\mu) - 2 \alpha \log |r|,
\end{align*}
so the $\mu$ that minimizes $E_R(\kap_*\mu)$ is the same as the $\mu$ that minimizes $E_Q(\mu)$.
\qedhere
\end{enumerate}
\end{proof}

\subsection{The extreme measure}
\label{S:extreme measure}

Let $c$ and $\Sigma_c$ be as in Proposition~\ref{P:potential version}.
Let $\mu_c$ be the equilibrium measure on $\Sigma_c$.
Lemma~\ref{L:evaluation of potential of mu_c} below shows that $\mu_c$ satisfies the requirements of Proposition~\ref{P:potential version}\eqref{I:extreme measures}.

\begin{lemma}
\label{L:evaluation of potential of mu_c}
\hfill
\begin{enumerate}[\upshape (a)]
\item \label{I:at 0} 
We have $U^{\mu_c}(0)=0$.
\item \label{I:value on Sigma_c} 
The function $U^{\mu_c}(z)$ is 
$\begin{cases} 
    -\log c &\textup{ for } z \in \Sigma_c, \\ 
    \le -\log c &\textup{ for } z \not \in \Sigma_c.
\end{cases}$
\item \label{I:real values} 
For all $r \in (-\infty,1]$, we have $U^{\mu_c}(r) = -\log M(r)$.
\end{enumerate}
\end{lemma}

\begin{proof}\hfill
\begin{enumerate}[\upshape (a)]
\item This holds for any measure supported on the unit circle, by definition of the potential.

\item The energy of the equilibrium measure on the circular arc $\Sigma_c$ is $-\log c$, as follows from \cite[Table~5.1]{Ransford95}. The inequality outside $\Sigma_c$ follows from nonnegativity of the Green function \cite[I.4]{Saff-Totik-1997}.  The equality on $\Sigma_c$ follows from the fact that the points on $\Sigma_c$ are regular points for the Dirichlet problem on $\C \setminus \Sigma_c$, as can be checked from \cite[Theorem~I.4.6]{Saff-Totik-1997}.

\item
By similar right triangles, the real part of either endpoint of $\Sigma_c$ is at distance $2c^2$ from $1$,
so $\kap(\Sigma_c) = 2[1-2c^2,1] = [2-4c^2,2]$.
Let $\ell(z) \colonequals c^2 z + 2-2c^2$, so $\ell([-2,2]) = [2-4c^2,2]$.
Let $\mu_\Sigma$ be the uniform probability measure on $\Sigma$,
so $U^{\mu_{\Sigma}}(z)=0$ on $\C_{\le 1}$ \cite[Example~0.5.7]{Saff-Totik-1997}.
By Lemma~\ref{L:pushforward}\eqref{I:compare equilibrium measures}, $\kap_*\mu_c$ and $\kap_*\mu_\Sigma$ are the equilibrium measures on $[2-4c^2,2]$ and $[-2,2]$, respectively, so $\kap_*\mu_c = \ell_* \kap_* \mu_\Sigma$.
Given $r \in (0,1]$, define $r' \in (0,1]$ by $\kap(r) = \ell(\kap(r'))$.  
Applying $U^{\kap_*\mu_c} = U^{\ell_* \kap_* \mu_\Sigma}$ yields
\begin{align*}
    U^{\kap_*\mu_c}(\kap(r)) &= U^{\kap_* \mu_\Sigma}(\kap(r')) - \log c^2 &&\textup{(since $\ell$ scales distances by $c^2$)} \\
    2 \, U^{\mu_c}(r) + \log r &= (2 \, U^{\mu_{\Sigma}}(r') + \log r') - \log c^2 &&\textup{(by Lemma~\ref{L:pushforward}\eqref{I:compare potentials} twice)} \\
    U^{\mu_c}(r) &= \tfrac12 \log (r'/r) - \log c &&\textup{(since $U^{\mu_{\Sigma}}(z)=0$ on $\C_{\le 1}$)} \\
    &= - \log M(r) &&\!\!\!\!\!\!\!\!\!\!\!\!\!\!\!\!\!\!\!\!\!\!\!\!\!\!\!\!\!\!\textup{(algebraic computation yields $r' \, M(r)^2 = r c^2$).}
\end{align*}
To extend to $(-\infty,1]$, observe that $U^{\mu_c}(r)$ and $-\log M(r)$ are real analytic on $(-\infty,1)$.
\qedhere
\end{enumerate}
\end{proof}

\subsection{Proof of optimality}

The idea for proving inequality \eqref{E:potential upper bound} is that it should be a nonnegative linear combination (really an integral) of the inequalities \eqref{E:c_assumption}.  The ``coefficients'' of the linear combination are given by a measure $\nu_\alpha$ belonging to a family that we describe now.
The $r=0$ and $r=1$ cases of \eqref{E:potential upper bound} follow from $M(0)=1$ and $M(1)=c$, so we assume $r \in (0,1)$.
For $\alpha \in \R_{\ge 0}$, let $\nu_\alpha$ be the equilibrium measure on $\Sigma$ for weight $Q_{\alpha}(z) \colonequals \alpha \log|z - r|$.

\begin{lemma}
\label{L:nu measures}
\hfill
\begin{enumerate}[\upshape (a)]
\item \label{I:theta(alpha)}
For every $\alpha \ge 0$, there exists $\theta(\alpha) \in (0, \pi] $ such that $\supp(\nu_{\alpha})$ is the arc $\{e^{it} \colon |t| \le \theta(\alpha)\}$.
\item \label{I:theta decreasing}
The function $\theta$ is decreasing and continuous.
Also, $\theta(0) = \pi$ and $\lim_{\alpha \to \infty} \theta(\alpha) = 0$.
\item \label{I:C}
Let $\alpha$ be such that $\supp(\nu_\alpha) = \Sigma_c$.
Then there is a constant $C$ such that 
\begin{equation}
\label{E:no_better}
   U^{\nu_\alpha}(z) + Q_\alpha(z)  \;\textup{ is }\; 
   \begin{cases}     C &\textup{ for all } z \in \Sigma_c, \\
                 \ge C &\textup{ for all } z \in \Sigma  \setminus \Sigma_c.
   \end{cases}
\end{equation}
\end{enumerate}
\end{lemma}

\begin{proof}
The literature contains similar results for an interval; to use them, we push forward by $\kap$.
By Lemma~\ref{L:pushforward}\eqref{I:compare equilibrium measures},
$\kap_*\nu_\alpha$ is the equilibrium measure on $[-2,2]$ for weight $Q_1(z) \colonequals \alpha \log|z-\kap(r)|$.
\begin{enumerate}[\upshape (a)]
\item 
We need to prove that $\supp(\kap_* \nu_\alpha) = [2\cos(\theta(\alpha)),2]$ for some $\theta(\alpha) \in (0,\pi]$.
Pushing forward $\kap_* \nu_\alpha$ by $z \mapsto 2-z$ gives the equilibrium measure on $[0,4]$ for weight 
$Q_2(z) \colonequals \alpha\log|\kap(r) - 2 + z|$. 
The function 
\[x \, Q_2'(x) = \frac{\alpha x}{\kap(r) - 2 + x}\]
is increasing on $[0, 4]$, 
so \cite[Theorem~IV.1.10(c)]{Saff-Totik-1997} implies that the support is an interval.
The corresponding interval for $\kap_* \nu_\alpha$ is contained in $[-2,2]$, and must contain $2$ since otherwise 
we could translate the measure right to reduce the weighted energy.
\item
By \cite[Theorem~IV.1.6(f)]{Saff-Totik-1997}, $\supp(\nu_\alpha)$ is decreasing and continuous,\footnote{Although \cite[Theorem~IV.1.6(f)]{Saff-Totik-1997} claims only right continuity, it can be applied with $Q$ replaced by $-Q$ to get left continuity.}
so $\theta$ is too.
Since $\nu_0$ is the uniform measure on $\Sigma$, we have $\theta(0)=\pi$.
For fixed $\epsilon>\epsilon'>0$, if $\kap_*\nu_\alpha$ has any mass to the left of $2-\epsilon$,
redistributing it according to the equilibrium measure on $[2-\epsilon',2]$
increases the energy by $O(1)$ but decreases the contribution from the weight 
by at least a positive constant times $\alpha$, 
so if $\alpha$ is sufficiently large, $\kap_*\nu_\alpha$ cannot have such mass; 
in other words, $\supp(\kap_*\nu_\alpha) \subset [2-\epsilon,2]$ for sufficiently large $\alpha$.
This holds for every $\epsilon$, so $\lim_{\alpha \to \infty} \theta(\alpha) = 0$.
\item
The number $\alpha$ exists by~\eqref{I:theta decreasing}.
Let $C$ be the modified Robin constant for $Q_\alpha$ \cite[p.~27]{Saff-Totik-1997}.
By \cite[Theorem~I.1.3(d,f)]{Saff-Totik-1997}, \eqref{E:no_better} holds outside a zero-capacity subset of $\Sigma$.
On the other hand, the points in $\Sigma$ are regular points for the Dirichlet problem in $\overline{\C} \setminus \Sigma$ by Wiener's theorem \cite[Theorem~I.4.6]{Saff-Totik-1997}, so \cite[Theorem~I.5.1(iv')]{Saff-Totik-1997} implies that $U^{\nu_\alpha}$ is continuous on $\Sigma$, as is $Q_\alpha$.
Thus \eqref{E:no_better} holds on all of $\Sigma$.\qedhere
\end{enumerate}
\end{proof}

\begin{proof}[Proof of \eqref{E:potential upper bound}]
By Lemma~\ref{L:Balay}, we may assume that $\mu$ is supported on the unit circle $\Sigma$, so $U^\mu(0)=0$.
Given $r \in (0,1)$, let $\alpha$ be as in Lemma~\ref{L:nu measures}\eqref{I:C}.
Then
\begin{align*}
-U^{\mu}(r) &= \frac{1}{\alpha} \int_{\Sigma}Q_{\alpha}(z) \, d\mu(z) &&\textup{(since $Q_{\alpha}(z) = \alpha \log |z-r|$)} \\
&\ge \frac{1}{\alpha} \left(C - \int_{\Sigma} U^{\nu_\alpha}(z) \, d\mu(z) \right) &&\textup{(by the inequality in~\eqref{E:no_better})} \\
&= \frac{C}{\alpha}  + \frac{1}{\alpha}\int_{\Sigma_c} \int_{\Sigma} \log|z - w|\, d\mu(z) \, d\nu_\alpha(w) && \textup{(by definition of $U^{\nu_\alpha}$)} \\
& = \frac{C}{\alpha}  - \frac{1}{\alpha}\int_{\Sigma_c} U^{\mu}(w) \, d\nu_\alpha(w) && \textup{(by definition of $U^{\mu}$)} \\
&\ge \frac {C}{\alpha} + \frac{1}{\alpha} \log c && \textup{(by \eqref{E:c_assumption} with $U^\mu(0)=0$).}
\end{align*}
When $\mu$ is $\mu_c$, Lemmas \ref{L:nu measures}\eqref{I:C} and~\ref{L:evaluation of potential of mu_c}\eqref{I:value on Sigma_c} show that both inequalities in this sequence are sharp, so 
\[
    -U^{\mu_c}(r) = \frac {C}{\alpha} + \frac{1}{\alpha} \log c \le -U^{\mu}(r).
\]
Thus $U^{\mu}(r) - U^{\mu}(0) = U^{\mu}(r) \le U^{\mu_c}(r) = - \log M(r)$, by Lemma~\ref{L:evaluation of potential of mu_c}\eqref{I:real values}.
\end{proof}

\begin{proof}[Proof of \eqref{E:potential lower bound}]
For $\beta \in \R_{\ge 0}$, let $\nu'_\beta$ be the equilibrium measure on $\Sigma$ for weight $R_{\beta}(z) \colonequals -\beta \log|z + r|$. As in Lemma \ref{L:nu measures}, there exists $\beta>0$ and a real constant $D$ such that $\supp(\nu'_{\beta}) = \Sigma_c$ and  
\[
    U^{\nu'_{\beta}}(z) + R_{\beta}(z)\;\textup{ is }\; 
   \begin{cases}    D  &\textup{ for all } z \in \Sigma_c, \\
                 \ge D &\textup{ for all } z \in \Sigma  \setminus \Sigma_c. 
    \end{cases}
\]
We may replace $\mu$ by the $\widehat{\mu}$ given by Lemma~\ref{L:Balay}, which implies that $U^\mu(r) - U^\mu(0) \ge U^{\widehat{\mu}}(r)$ for every $r \in [0,\infty)$.
The rest of the proof is entirely analogous to the proof of \eqref{E:potential upper bound}.
\end{proof}

%****************************************************************************
\section*{Acknowledgments} 

We thank Francesc Fit\'e, Everett Howe, and Stefano Marseglia for discussions.
We also thank Everett Howe and Felipe Voloch for suggesting the Weil references in the first sentence.

%\bibliographystyle{alpha}
%\bibliography{biblio}
\printbibliography
\end{document}